\documentclass{article}
\usepackage{amssymb}
\usepackage{amsmath}
\usepackage{amsfonts}
\usepackage{latexsym}
\usepackage{epsf,graphicx,subfigure}

\newtheorem{theorem}{Theorem}

\newtheorem{corollary}[theorem]{Corollary}

\newtheorem{definition}[theorem]{Definition}
\newtheorem{example}[theorem]{Example}

\newtheorem{lemma}[theorem]{Lemma}

\newtheorem{proposition}[theorem]{Proposition}
\newtheorem{remark}[theorem]{Remark}

\newenvironment{proof}[1][Proof]{\noindent\textbf{#1.} }{\ \rule{0.5em}{0.5em}}
\newenvironment{proof7}[1][Proof of Theorem 7]{\noindent\textbf{#1.} }{\ \rule{0.5em}{0.5em}}
\newenvironment{proof12}[1][Proof of Theorem 12]{\noindent\textbf{#1.} }{\ \rule{0.5em}{0.5em}}
\newenvironment{proof18}[1][Proof of Theorem 18]{\noindent\textbf{#1.} }{\ \rule{0.5em}{0.5em}}
\newenvironment{proof39}[1][Proof of Theorem 39]{\noindent\textbf{#1.} }{\ \rule{0.5em}{0.5em}}
\newcommand {\Hol}{\mathop{\rm Hol}\nolimits}
\renewcommand{\Re}{\mathop{\rm Re}\nolimits}
\renewcommand{\Im}{\mathop{\rm Im}\nolimits}

\begin{document}

\title{Linearization models for parabolic dynamical systems via Abel's
functional equation}
\author{Mark Elin
\\ {\small Department of Mathematics}
\\ {\small ORT  Braude College}
\\ {\small P.O. Box 78, 21982 Karmiel, Israel}
\\ {\small e-mail: mark$\_$elin@braude.ac.il}
\\ Dmitry Khavinson
\\ {\small  Department of Mathematics}
\\ {\small  University of South Florida}
\\ {\small  Tampa, FL 33620-5700, USA}
\\ {\small e-mail: dkhavins@cas.usf.edu}
\\ Simeon Reich
\\ {\small Department of Mathematics}
\\ {\small The Technion --- Israel Institute of Technology}
\\ {\small 32000 Haifa, Israel}
\\ {\small e-mail: sreich@tx.technion.ac.il}
\\ David Shoikhet
\\ {\small Department of Mathematics}
\\ {\small ORT Braude College}
\\ {\small P.O. Box 78, 21982 Karmiel, Israel}
\\ {\small e-mail: davs@braude.ac.il}}
\date{ }
\maketitle

\begin{abstract}
We study linearization models for continuous one-parameter
semigroups of parabolic type. In particular, we introduce new
limit schemes to obtain solutions of Abel's functional equation
and to study asymptotic behavior of such semigroups. The crucial
point is that these solutions are univalent functions convex in
one direction. In a parallel direction, we find analytic
conditions which determine certain geometric properties of those
functions, such as the location of their images  in either a
half-plane or a strip, and their containing either a half-plane or
a strip. In the context of semigroup theory these geometric
questions may be interpreted as follows: is a given one-parameter
continuous semigroup either an outer or an inner conjugate of a
group of automorphisms? In other words, the problem is finding a
fractional linear model of the semigroup which is defined by a
group of automorphisms of the open unit disk. Our results enable
us to establish some new important analytic and geometric
characteristics of the asymptotic behavior of one-parameter
continuous semigroups of holomorphic mappings, as well as to study
the problem of existence of a backward flow invariant domain and
its geometry.
\end{abstract}

\section{\protect\bigskip Introduction}

Let $\Delta $ be the open unit disk in the complex plane
$\mathbb{C}$ and let $\Omega $ be a subset of $\mathbb{C}$. By
$\Hol (\Delta ,\Omega )$ we denote the set of all holomorphic
functions (mappings) from $\Delta$ into $\Omega$.

Linearization models for continuous semigroups of holomorphic
mappings in various settings have an extensive history, beginning
with the study of continuous stochastic Markov branching processes
(see, for example, \cite{Harris}), Kolmogorov's backward equation
in probability theory and functional differential equations.

The linearization models for one-parameter semigroups have also
proved to be very useful in the study of composition operators and
their spectra (see, for example, \cite{BE-PH}, \cite{Cowen},
\cite{Siskakis}, \cite{CCC-MBD}). In the years since these classic
works appeared, much more information has been found.

A deep investigation of the behavior of one-parameter semigroups
near their boundary Denjoy-Wolff point based on geometric
properties of the model obtained by Abel's functional equation has
recently been conducted in \cite{CDP, CM-DS-CP1}.

It turns out that in those settings solutions of Abel's functional
equation are univalent functions convex in one-direction.

In a parallel direction, the class of functions convex in one
direction has been studied by many mathematicians (see, for
example, \cite{H-S}, \cite{Ci1}, \cite{Ci2}, \cite{GAW} and
\cite{Le}) as a subclass of the class of functions introduced by
Robertson in \cite{Rob1}.

Following this point of view, we study in the present paper
several additional geometric and analytic properties of functions
convex in one direction by using the asymptotic behavior of one
parameter semigroups near their attractive (or repelling) boundary
fixed points.

In an opposite direction this study yields new information
concerning some special (but very wide and important) classes of
complex dynamical systems generated by holomorphic functions
having fractional derivatives at their Denjoy-Wolff points.

\begin{definition} \label{def1}
A univalent function $h\in \Hol(\Delta ,\mathbb{C} )$ is said to
be convex in the positive direction of the real axis if for each
$z\in \Delta $ and $t\geq 0$, the point $h(z)+t$ belongs to
$h(\Delta )$.
\end{definition}

It is well known (see, for example, \cite{Shapiro}, p. 162) that
for each $z\in \Delta $, the limit
\begin{equation} \label{N1}
\lim_{t\rightarrow \infty }h^{-1}(h(z)+t)=:\zeta
\end{equation}
exists and belongs to $\partial \Delta $.

Moreover, since the family $S=\{F_{t}\}_{t\geq 0}$ defined by
\begin{equation} \label{N2}
F_{t}(z)=h^{-1}(h(z)+t)
\end{equation}
forms a one-parameter continuous semigroup of holomorphic
self-mappings of $\Delta $, it follows from the continuous version
of the Denjoy-Wolff Theorem (see \cite{BE-PH}, \cite{RS-SD-97b}
and \cite{SD}) that the limit point $\zeta $ in (\ref{N1}) is
unique and does not depend on $z\in \Delta $.

Without loss of generality we may set $\zeta =1$.

We denote by $\Sigma [1]$ the class of functions convex in the
positive direction of the real axis, normalized by the conditions
\begin{equation} \label{N3}
\lim_{t\rightarrow \infty }h^{-1}(h(z)+t)=1\text{ \ and \ }h(0)=0.
\end{equation}

In the reverse direction one can assign to each semigroup $
S=\{F_{t}\}_{t\geq 0}$ of holomorphic self-mappings of $\Delta $
with a boundary Denjoy-Wolff point $\zeta =1$, a univalent
function $h\in \Hol(\Delta , \mathbb{C})$ which is a common
solution to Abel's functional equations
\begin{equation} \label{V1}
h(F_{t}(z))=h(z)+t\text{, \ \ }z\in \Delta \text{, }t\geq 0\text{,}
\end{equation}
and hence is convex in the positive direction of the real axis
(see the proof of Proposition 3 below). The set $h\left( \Delta
\right)$ is called a \textit{planar domain} for $S$ and the pair
$(h,h(\Delta))$ is said to be a \textit{linearization model for
}$S$.

To be more precise, we recall that by the Berkson-Porta theorem
\cite{BE-PH} (see also \cite{AM-92}, \cite{RS-SD-96},
\cite{RS-SD-97} and \cite{RS}), for each one-parameter continuous
semigroup $S=\{F_{t}\}_{t\geq 0}$ of holomorphic self-mappings of
$\Delta$, the limit
\begin{equation}
f(z):=\lim_{t\rightarrow 0^{+}}\frac{z-F_{t}(z)}{t}
\end{equation}
exists. This limit function is called the \textit{infinitesimal
generator of} $S$. Moreover, the semigroup $S=\{F_{t}\}_{t\geq 0}$
can be reproduced as the solution of the Cauchy problem
\begin{equation}
\left\{
\begin{array}{c}
\dfrac{\partial u(t,z)}{\partial t}+f(u(t,z))=0 \\
\\
u(0,z)=z,
\end{array}
\right.
\end{equation}
where we set $F_{t}(z):=u(t,z)$, $t\geq 0$, $z\in \Delta$.

By $\mathcal{G}[1]$ we denote the class of functions which
consists of all the holomorphic generators $f$ of continuous
semigroups having the Denjoy-Wolff point $\zeta =1$. In this case
$f$ admits the Berkson-Porta representation
\begin{equation*}
f(z)=-(1-z)^{2}p(z)\text{, \ }z\in \Delta \text{,}
\end{equation*}
with $\Re p(z)\geq 0$ everywhere (see \cite{BE-PH} and \cite{SD}).

In addition, if $S=\{F_{t}\}_{t\geq 0}$ is the semigroup generated
by $f$ ($S$ is defined via the Cauchy problem (\ref{V3'}), where
we set $F_{t}(z):=u(t,z)$, $t\geq 0$, $z\in \Delta$), then
$\lim_{t\rightarrow \infty }F_{t}(z)=1$ and for each $t\geq 0$,
\begin{equation*}
\angle \lim_{z\rightarrow 1}\frac{\partial F_{t}(z)}{\partial z}=e^{-t\beta
},
\end{equation*}
where $\beta :=\angle \lim_{z\rightarrow 1}\dfrac{f(z)}{z-1}=\angle
\lim_{z\rightarrow 1}f^{\prime }(z)\geq 0$ (see, for example, \cite{E-S1}
and \cite{CDP}).

If $\beta >0$, then $S$ and its generator $f$ are said to be of
\textit{hyperbolic type}. Otherwise, that is, if $\beta =0$, the
semigroup $S$ and its generator $f$ are said to be of
\textit{parabolic type}.

Conversely, it is shown in \cite{E-S1} that if $f$ is a
holomorphic generator on $\Delta $ such that the angular limit
\begin{equation} \label{V4}
\angle \lim_{z\rightarrow 1}\frac{f(z)}{z-1}=:f^{\prime }(1)
\end{equation}
exists finitely with $\Re f^{\prime }(1)\geq 0$, then $f^{\prime
}(1)$ is, in fact, a real number, $f$  has no null points in
$\Delta $ and belongs to $\mathcal{G}[1]$.

This, in turn, implies (see Section 1 for details) that the formula
\begin{equation*}
h(z)=-\int\limits_{0}^{z}\frac{dz}{f(z)}
\end{equation*}
establishes a one-to-one correspondence between the classes $\mathcal{G}[1]$
and $\Sigma \lbrack 1]$.

In this paper we are interested in finding simple analytic conditions which
determine certain geometric properties of functions $\ h$ in the class $%
\Sigma [1]$, such as the location of the image $h(\Delta )$ in
either a half-plane or a strip, and its containing either a half
plane or a strip.

In the context of semigroup theory these geometric questions may be
interpreted as follows:

Is a given one-parameter continuous semigroup either an outer or
an inner conjugate of a group of automorphisms (see the definition
below)? In other words, the problem is finding a fractional linear
model of the semigroup which is defined by a group of
automorphisms of $\Delta $.

To be more precise, we need the following definition.

\begin{definition} \label{defV}
Let $S=\{F_{t}\}_{t\geq 0}$ be a one-parameter continuous
semigroup of holomorphic self-mappings of $\Delta $.

$\bullet $ A univalent mapping $\psi :\Delta \rightarrow \Delta $
is called an outer conjugator of $S$ if there is a one-parameter
semigroup $\{G_{t}\}_{t\geq 0}\subset \Hol(\Delta ,\Delta )$ of
linear fractional transformations of $\Delta $ such that
\begin{equation}
\psi \circ F_{t}=G_{t}\circ \psi \text{.}  \label{V2}
\end{equation}

$\bullet $ A univalent mapping $\varphi :\Delta \rightarrow \Delta
$ is called an inner conjugator of $S$ if there is a one-parameter
continuous semigroup $\{G_{t}\}_{t\geq 0}\subset \Hol(\Delta
,\Delta )$ of linear fractional transformations of $\Delta $ such
that
\begin{equation} \label{V3'}
F_{t}\circ \varphi =\varphi \circ G_{t}\text{.}
\end{equation}
\end{definition}

Note that by the Linear Fractional Model theorem (see \cite{BP-SJ}
and \cite{Shapiro}) for each holomorphic self-mapping there exists
a conjugating function $G:\Delta \mapsto \mathbb{C} $ which is not
necessarily a self-mapping of $\Delta $.

However, in our considerations the main point in Definition 2 is that the
conjugators $\varphi $ and $\psi $ are \textit{self}-mappings of $\Delta $.

It is rather a transparent fact (see Sections 2 and 3 for details)
that if $S=\{F_{t}\}_{t\geq 0}$ is a semigroup with a boundary
Denjoy-Wolff point $\zeta \in \partial \Delta $ (say $\zeta =1$)
and $h$ is the common solution to Abel's equations (\ref{V1})
normalized by (\ref{V3'}), then an outer (respectively, inner)
conjugator of $S$ exists if and only if $h(\Delta )$ lies in
(respectively, contains) a half-plane $\Pi $.

In both cases the conjugate semigroup $\{G_{t}\}_{t\geq 0}$ of
linear fractional transformations (see formulas (\ref{V2}) and
(\ref{V3'})) is a group of automorphisms of $\Delta $ if and only
if $\Pi $ is horizontal, i.e., its boundary is parallel to the
real axis.

The latter situation is of special interest for inner conjugation
since equation (\ref{V3'}) means that the semigroup
$S=\{F_{t}\}_{t\geq 0}$ forms a group of automorphisms of the
domain $\varphi (\Delta )$ which is called a \textit{backward flow
invariant domain} for $S$ (see \cite{E-S-Z1}). We will concentrate
on this issue in Section 3 for semigroup generators which are
fractionally differentiable at their boundary Denjoy-Wolff points.
In the meantime we need the following observation.

It can be shown (see \cite{E-D-2006} and Lemma 3 below) that if
$h$ belongs to $\Sigma [1]$, then the angular limit
\begin{equation} \label{ang}
\mu :=\angle \lim_{z\rightarrow 1}(1-z)h^{\prime }(z)
\end{equation}
exists and belongs to $\left( 0,\infty \right] $.

It turns out that this limit $\mu $ is finite if and only if
$\left\vert \Im h\left( z\right) \right\vert $ is bounded.
Moreover, the number $\pi \mu $ is the width of the minimal
horizontal strip which contains $h(\Delta )$.

If $\mu $ in (\ref{ang}) is infinite, it may happen that for some
$\alpha >0$, the angular limit
\begin{equation} \label{N11}
\angle \lim_{z\rightarrow 1}(1-z)^{1+\alpha }h^{\prime }(z)=:\mu ,
\end{equation}%
or even the unrestricted limit
\begin{equation} \label{unr}
\lim_{\substack{ z\rightarrow 1  \\ z\in \Delta }}(1-z)^{1+\alpha
}h^{\prime }(z)=:\mu ,
\end{equation}
exist finitely.

The latter condition, for example, follows from the geometric
property of $h(\Delta )$ to have a Dini smooth corner of opening
$\pi \alpha $ at infinity for some $\alpha \in \left( 0,2\right] $
(see \cite{PC-92}).

Therefore it is natural to consider the subclasses $\Sigma
_{A}^{\alpha }[1]$ , respectively, $\Sigma ^{\alpha }[1]$ ($\Sigma
^{\alpha }[1]\subset \Sigma _{A}^{\alpha }[1]$), $\alpha \geq 0$,
which consist of those $h\in \Sigma [1]$ satisfying (\ref{N11}),
respectively, (\ref{unr}).

We will first show that if $\Sigma _{A}^{\alpha }[1]\neq
\varnothing $, then $\alpha $ must belong to the interval $\left[
0,2\right] $. Furthermore, if $\alpha >0$ and $h\in \Sigma
^{\alpha }[1]$, then the boundedness from above (or from below) of
$\Im h\left( z\right) $ implies some necessary simple relations
between the numbers $\alpha $ and $\mu $ in (\ref{unr}). Note that
if $\alpha >0$, then $\Im h\left( z\right) $ cannot be bounded
from both above and below.

Finally, we would like to emphasize that these relations, as well
as the geometric properties of these classes, are connected to the
asymptotic behavior of the semigroups defined by (\ref{N2}). Also,
some special situations arise when $\alpha $ attains its integer
values $\{0,1,2\}$.

\section{Subordination theorems (outer conjugation)}

\subsection{Auxiliaries results}

 Analytic descriptions of functions convex in one direction have
been presented by many mathematicians using different methods of
classical geometric function theory including the prime ends
theory and Julia's lemma (see, for example, \cite{Rob1},
\cite{H-S}, \cite{Ci1}, \cite{Ci2}, \cite{GAW} and \cite{Le}). For
completeness we begin with the following description of the class
$\Sigma [1]$, which is a slight generalization of those given in
\cite{H-S}, \cite{Ci1}, \cite{Ci2}, and \cite{Le}, and is based on
the relations between this class and the class $\mathcal{G}[1]$.

\begin{proposition} \label{lem1}
A function $h\in \Hol(\Delta , \mathbb{C} )$, $h\neq 0$, with
$h(0)=0$, belongs to $\Sigma [1]$ if and only if
\begin{equation}  \label{N4}
\Re \left[ (1-z)^{2}h^{\prime }(z)\right] \geq 0
\end{equation}
and if and only if the function $\ f\in \Hol(\Delta ,\Omega )$,
defined by $f(z)=-\dfrac{1 }{h^{\prime }(z)}$, belongs to
$\mathcal{G}[1]$.
\end{proposition}

\begin{proof}
\label{proofN8} Indeed, if $h\in \Sigma [1]$, then $
S=\{F_{t}\}_{t\geq 0}$, where $F_{t}$ is defined by (\ref{N2}), is
a one-parameter continuous semigroup with the Denjoy-Wolff point
$\zeta =1$. Therefore, by the Berkson-Porta theorem \cite{BE-PH},
its infinitesimal generator
\begin{equation} \label{N5}
f(z):=\lim_{t\rightarrow 0^{+}}\frac{z-F_{t}(z)}{t}
\end{equation}
admits the representation
\begin{equation}
f(z)=-(1-z)^{2}p(z)  \label{N6}
\end{equation}
with
\begin{equation} \label{N7}
\Re p(z)\geq 0\text{.}
\end{equation}

On the other hand, differentiating (\ref{N2}) at $t=0^{+}$, we \ obtain
\begin{equation*}
f(z)=-\frac{1}{h^{\prime }(z)}\text{.}
\end{equation*}%
Thus (\ref{N6}) and (\ref{N7}) imply (\ref{N4}).

Conversely, assume now that $h\in \Hol(\Delta ,\mathbb{C} )$,
$h\neq 0$, with $h(0)=0$, satisfies (\ref{N4}). If $\Re \left[
(1-z)^{2}h^{\prime }(z)\right] =0$, then $h(z)=\dfrac{ibz}{1-z}$
for some $b\in \mathbb{R} $, and clearly $h$ belongs to $\Sigma
[1]$. In this case the function $\ f\in \Hol(\Delta ,\Omega )$,
defined by $f(z)=-\dfrac{1}{h^{\prime }(z)},$ is a generator of a
group of parabolic automorphisms of $\Delta$, hence belongs to
$\mathcal{G}[1]$.

So we may assume that
\begin{equation*}
\Re \left[ (1-z)^{2}h^{\prime }(z)\right] >0.
\end{equation*}
In this case, using the convex function $g(z)=\frac{z}{1-z}$, we
have
\begin{equation*}
\Re \left[ \frac{h^{\prime }(z)}{g^{\prime }(z)}\right] >0,
\end{equation*}
which means that $h$ is close-to-convex. Hence it must be univalent (see
\cite{DP}, Theorem 2.17, p. 47).

In addition, we know that the function $f\in \Hol(\Delta ,
\mathbb{C} )$ defined by
\begin{equation} \label{N8}
f(z)=-\frac{1}{h^{\prime }(z)}
\end{equation}
admits the Berkson-Porta representation (\ref{N6}) with
(\ref{N7}). Therefore, the Cauchy problem
\begin{equation}  \label{N9}
\left\{
\begin{array}{c}
\dfrac{\partial u(t,z)}{\partial t}+f(u(t,z))=0 \\
\\
u(0,z)=z
\end{array}
\right.
\end{equation}
is solvable with $\left\vert u(t,z)\right\vert <1$ for each $z\in
\Delta $ and all $t\geq 0$. Moreover, it follows from (\ref{N8})
and (\ref{N9}) that
\begin{equation*}
h^{\prime }(u(t,z))\cdot \frac{\partial u(t,z)}{\partial t}=1
\end{equation*}
for all $z\in \Delta $ and $t\geq 0$.

Integrating this equality, we get
\begin{equation*}
h(u(t,z))=h(z)+t\in h(\Delta )
\end{equation*}
for all $z\in \Delta $ and $t\geq 0$. Hence $h$ is indeed convex in the
positive direction of the real axis.

Finally, again by the Berkson-Porta theorem, we have
$\lim_{t\rightarrow \infty }u(t,z)=1$. Hence $h\in \Sigma [1]$ and
the proof is compete.
\end{proof}

Thus, we have shown that the formula
\begin{equation*}
h(z)=-\int\limits_{0}^{z}\frac{dz}{f(z)}
\end{equation*}
determines a one-to-one correspondence between the classes
$\mathcal{G}[1]$ and $\Sigma [1]$.

For $h\in \Sigma [1]$ we assume now that $h(\Delta )$ lies in a
horizontal half-plane $\Pi =h_{1}(\Delta )$, where $h_{1}:\Delta
\mapsto \mathbb{C} $ is defined by
\begin{equation*}
h_{1}(z)=\frac{ibz}{1-z}
\end{equation*}
for some $b\in \mathbb{R} $.

Consider the mapping $\psi :\Delta \mapsto \Delta $ defined by
\begin{equation}
\psi (z)=h_{1}^{-1}(h(z)),\text{ \ \ }\psi (0)=0\text{,}  \label{n2''}
\end{equation}

or explicitly
\begin{equation} \label{n3''}
\psi (z)=\frac{h(z)}{ib+h(z)}\text{.}
\end{equation}
Define a group of (parabolic) automorphisms of $\Delta $ by
\begin{eqnarray}
G_{t}(z) &=&h_{1}^{-1}(h_{1}(z)+t),\text{ \ \ }t\in \mathbb{R},
\label{n4''} \\ &=&\frac{h_{1}(z)+t}{ib+h_{1}(z)+t}=  \notag \\
&=&\frac{\dfrac{ibz}{1-z}+t}{ib+\dfrac{ibz}{1-z}+t}=\frac{ibz+t(1-z)}{
ib+t(1-z)}\text{.}  \notag
\end{eqnarray}
It is clear that
\begin{equation*}
\lim_{t\rightarrow \pm \infty }G_{t}(z)=1\text{.}
\end{equation*}

Now for $t\geq 0$ we have by (\ref{n4''})
\begin{eqnarray*}
G_{t}(\psi (z)) &=&h_{1}^{-1}(h_{1}(\psi (z))+t)= \\
&=&h_{1}^{-1}(h(z)+t)=h_{1}^{-1}(h(F_{t}(z)))=\psi (F_{t}(z))\text{,}
\end{eqnarray*}
where $F_{t}(z)=h^{-1}(h(z)+t)$.

Thus, for $S=\{F_{t}\}_{t\geq 0}$ generated by $f$, defined by
(\ref{N8}), there is a group $\{G_{t}\}_{t\in \mathbb{R}}$ of
parabolic automorphisms of $\Delta $, such that $S$ is an outer
conjugate of $\{G_{t}\}_{t\geq 0}$  by $\psi :\Delta \mapsto
\Delta $. Moreover, the domain $\Omega =\psi (\Delta )$ is
invariant for the semigroup $\{G_{t}\}_{t\geq 0}$.

Conversely, let $\{G_{t}\}_{t\geq 0}$ be a group of parabolic
automorphisms of $\Delta $ such that $\lim_{t\rightarrow \pm
\infty }G_{t}(z)=1$, and assume that there is a holomorphic
Riemann mapping $\psi $ of $\Delta $ into $\Delta $, $\psi (0)=0$,
such that
\begin{equation}  \label{n5''}
G_{t}(\psi (z))=\psi (F_{t}(z))
\end{equation}
for all $z\in \Delta $ and $t\geq 0$.

Differentiating (\ref{n5''}) with respect to $t\geq 0$, we have
\begin{equation} \label{n6''}
g(\psi (z))=\psi ^{\prime }(z)\cdot f(z),
\end{equation}
where $g$ and $f$ are generators of $\{G_{t}\}_{t\geq 0}$ and $
\{F_{t}\}_{t\geq 0}$, respectively.

Explicitly, $g(z)$ can be calculated as follows:
\begin{eqnarray}
g(z) &=&\lim_{t\rightarrow
0^{+}}\frac{z-G_{t}(z)}{t}=\lim_{t\rightarrow 0}
\frac{z-\dfrac{ibz+t(1-z)}{ib+t(1-z)}}{t}=  \label{n7''} \\
&=&\lim_{t\rightarrow
0}\frac{ibz+tz(1-z)-ibz-t(1-z)}{t(ib+t(1-z))}=  \notag
\\
&=&\lim_{t\rightarrow
0^{+}}\frac{t(1-z)(z-1)}{t(ib+t(1-z))}=-\frac{(1-z)^{2}
}{ib}\text{. \ }  \notag
\end{eqnarray}

Since $f(z)=-\dfrac{1}{h^{\prime }(z)}$, we get from (\ref{n6''})
and (\ref {n7''}) that
\begin{equation} \label{n8''}
h^{\prime }(z)=\frac{\psi ^{\prime }(z)ib}{(1-\psi (z))^{2}}\text{.}
\end{equation}

Integrating (\ref{n8''}) from $0$ to $z$, we obtain
\begin{equation} \label{n9''}
h(z)=\frac{ib\psi (z)}{1-\psi (z)}\text{.}
\end{equation}

Now, setting $w=\psi (z)\in \Omega $, we get from (\ref{n9''})
that
\begin{equation*}
h(\psi ^{-1}(w))=\frac{ibw}{1-w}=:h_{1}(w)\text{.}
\end{equation*}

Since this function can be extended to the whole disk $\Delta $
and $h_{1}(\Delta )$ is a horizontal half-plane, equality
(\ref{n9''}) means that $h(\Delta )\subset h_{1}(\Delta )$. Thus
we have proved the following assertion.

\begin{proposition} \label{Prop01}
Let $\{F_{t}\}_{t\geq 0}$ be a semigroup with the Denjoy-Wolff
point $z=1$, such that $F_{t}(z)=h^{-1}(h(z+t))$, $h\in \Sigma
[1]$. Then $h(\Delta )$ lies in a horizontal half-plane if and
only if there is a group $\{G_{t}\}$, $t\in \mathbb{R}$, of
parabolic automorphisms of $\Delta $ and a conformal self-mapping
$\psi $ of $\Delta $ such that for all $t\geq 0$, the semigroup
$\{F_{t}\}$ outer conjugates with $\{G_{t}\}$:
\begin{equation*}
\psi (F_{t}(z))=G_{t}(\psi (z)),\text{ \ \ }t\geq 0\text{.}
\end{equation*}
\end{proposition}

\begin{remark} \label{Remark01}
In particular, each hyperbolic group can conjugate with a
parabolic one. This can easily be shown if we translate our
consideration to the right half-plane.

Namely, let $\{\widehat{F}_{t}\}$ be the group of hyperbolic
automorphisms of $\Pi _{+}=\{w\in \mathbb{C} :\Re w>0\}$ of the
form
\begin{equation*}
\widehat{F}_{t}(w)=e^{t}w\text{.}
\end{equation*}

Consider the univalent function $\widehat{\psi }:\Pi _{+}\mapsto
\mathbb{C}$ given by
\begin{equation*}
\widehat{\psi }(w)=i\log w+r\text{,}
\end{equation*}%
where $r>-\dfrac{\pi }{2}$.

Then for $w\in \Pi _{+}=\left\{ -\frac{\pi }{2}<\arg w<\frac{\pi
}{2} \right\} $, we have $\ \dfrac{\pi }{2}+r>\Re \widehat{\psi
}(w)=\arg w+r>0$.

So, $\widehat{\psi }$ maps $\Pi _{+}$ into a vertical strip
$\widehat{\Omega }\subset \Pi _{+}$.

In addition,
\begin{eqnarray*}
\widehat{\psi }(\widehat{F}_{t}(z)) &=&i\log (e^{t}w)+r= \\
&=&it+i\log w+r=it+\widehat{\psi }(w)\text{.}
\end{eqnarray*}

Thus, if we define a group of parabolic automorphisms of $\Pi _{+}$ by
\begin{equation*}
\widehat{G}_{t}(w)=w+it\text{,}
\end{equation*}
we get that
\begin{equation*}
\widehat{\psi }(\widehat{F}_{t}(w))=\widehat{G}_{t}(\widehat{\psi
}(w))\text{ .}
\end{equation*}

Applying now the Cayley transform $C:\Delta \mapsto \Pi ^{+}$,
$C(z)=\dfrac{1+z}{1-z}$, and setting $F_{t}=C^{-1}\circ
\widehat{F}_{t}\circ C$, $G_{t}=C^{-1}\circ \widehat{G}_{t}\circ
C$ and $\psi =C^{-1}\circ \widehat{\psi }\circ C$, one obtains the
needed relation in the unit disk
\begin{equation*}
\psi (F_{t}(z))=G_{t}(\psi (z))\text{, \ \ \ \ }z\in \Delta \text{.}
\end{equation*}

In this case $\Omega =C^{-1}(\widehat{\Omega })$ is the domain bounded by
two horocycles internally tangent to $\partial \Delta $ at $z=1$.
\end{remark}

To proceed we need the following lemma.

\begin{lemma} \label{lem2}
Let $h\in \Sigma [1]$. Then the angular limit
\begin{equation*}
\angle \lim_{z\rightarrow 1}(1-z)h^{\prime }(z):=\mu
\end{equation*}
exists and is either a positive real number or infinity.
\end{lemma}

\begin{proof} \label{proofN7}
If we denote $p(z):=1/(1-z)^{2}h^{\prime }(z)$, then we have $\Re
p(z)\geq 0$.

Now it follows from the Riesz-Herglotz representation of $p$,
\begin{equation*}
p(z)=\oint\limits_{\partial \Delta }\frac{1+\overline{\zeta
}z}{1-\overline{ \zeta }z}dm(\zeta )+i\gamma ,
\end{equation*}
where $m$ is a positive measure on $\partial \Delta $ and $\gamma
\in \mathbb{R}$, that
\begin{equation*}
1/\mu =\angle \lim_{z\rightarrow 1}(1-z)p(z)=2m(1)\in \lbrack 0,\infty ).
\end{equation*}%
This proves our assertion.
\end{proof}

\subsection{Main results}

We begin with a description of the class $\Sigma _{A}^{0}[1].$

\begin{theorem} \label{teor1}
Let $h\in \Sigma [1]$. The following assertions are equivalent:

(i) $h\in \Sigma _{A}^{0}[1]$, i.e.,
\begin{equation*}
\angle \lim_{z\rightarrow 1}(1-z)h^{\prime }(z)=\mu
\end{equation*}
exists finitely;

(ii) $\Omega =h(\Delta )$ lies in a horizontal strip;

(iii) $h$ is a Bloch function, i.e.,
\begin{equation*}
\left\Vert h\right\Vert _{B}:=\sup_{z\in \Delta }(1-\left\vert z\right\vert
^{2})\left\vert h^{\prime }(z)\right\vert <\infty .
\end{equation*}

Moreover, in this case,

(a) the smallest strip which contains $h(\Delta )$ is
\begin{equation*}
\left\{ w\in \mathbb{C} :\left\vert \Im w-a\right\vert <\frac{\pi
\mu }{2}\right\} ,
\end{equation*}
where
\begin{equation*}
a=\lim_{r\rightarrow 1^{-}}\Im h(r).
\end{equation*}

(b) \ $2\left\vert \mu \right\vert \leq \left\Vert h\right\Vert _{B}\leq
4\left\vert \mu \right\vert $.
\end{theorem}

\begin{remark} \label{rem1}
One of the tools we use in the proof of our results is the
Koebe Distortion Theorem (see, for example, \cite{PC-92}) which asserts:

if $h$ is a univalent function on $\Delta $, then
\begin{equation} \label{R1}
\frac{1}{4}(1-\left\vert z\right\vert ^{2})\left\vert h^{\prime
}(z)\right\vert \leq \delta (h(z))\leq (1-\left\vert z\right\vert
^{2})\left\vert h^{\prime }(z)\right\vert ,\text{ \ }z\in \Delta \text{,}
\end{equation}
where
\begin{equation*}
\delta (w)=dist\left( w,\partial h(\Delta )\right) \text{, \ }w\in
h(\Delta ) \text{.}
\end{equation*}

Once the equivalence of (i)-(iii) and assertion (a) are proved, we
have by (\ref{R1}) the following estimate:
\begin{equation} \label{R2}
\frac{\pi \left\vert \mu \right\vert }{2}\leq \left\Vert
h\right\Vert _{B}\leq 2\pi \left\vert \mu \right\vert .
\end{equation}

Clearly, assertion (b) improves (\ref{R2}).

Such estimates are very important in finding non-trivial bounds for integral
means of $\left\vert h(z)\right\vert $, as well as its even powers (see, for
example, \cite{PC-92}, p. 186).
\end{remark}

\begin{proof7}
The equivalence of assertions (i) and (ii), as well as assertion
(a), are proved in \cite{E-D-2006}. The implication (ii)
$\Longrightarrow $(iii) follows immediately from (\ref{R1}).

Now we will show that (iii) implies (ii). Assume to the contrary
that there is no horizontal strip which contains $h(\Delta )$.
Then for each $M>0$, there are two points $w_{1}$ and $w_{2}$ in
$h(\Delta )$ such that
\begin{equation} \label{R3}
\Im w_{1}-\Im w_{2}>2M.
\end{equation}

Take any continuous curve $\Gamma \subset h(\Delta )$ ending at
the points $w_{1}$ and $w_{2}$, and consider the half-strip
\begin{equation*}
\Omega _{1}=\left\{ w\in \mathbb{C}:\Re w>b,\text{ \ }\Im
w_{2}<\Im w<\Im w_{1}\right\},
\end{equation*}%
where $b=\max_{w\in \Gamma }\Re w$.

Since $h\in \Sigma (1)$, this half-strip $\Omega _{1}$ must lie in $\Omega
=h(\Delta )$.

If now $l$ is the midline of $\Omega _{1}$, then one can find $w\in l$ such
that $\delta (w)=dist(w,\partial \Omega )>M$.

Since $M$ is arbitrary, this contradicts (iii) by (\ref{R1}). Thus (iii)
does indeed imply (ii).

Finally, it remains to be shown that assertion (b) holds. The
left-hand side of the inequality in (b) is immediate. To prove the
right-hand side of this inequality, we consider the function
$\widehat{h}\in \Hol(\Delta , \mathbb{C} )$ defined by
$\widehat{h}(z)=\exp \left( -\dfrac{1}{\mu }h(z)\right) $. It has
been shown in \cite{E-D-2006} that this function is a univalent
(starlike with respect to a boundary point) function on $\Delta $
with $h(z)\neq 0$, $z\in \Delta$. Therefore it follows from
Proposition 4.1 in \cite{PC-92} that $g(z):=\log \widehat{h}(z)$
is a Bloch function with
\begin{equation*}
\left\Vert g(z)\right\Vert _{B}\leq 4\text{.}
\end{equation*}
But $g(z)=-\dfrac{1}{\mu }h(z)$. Hence
\begin{equation*}
\left\Vert h(z)\right\Vert _{B}\leq 4\left\vert \mu \right\vert \text{.}
\end{equation*}
\end{proof7}

So, according to relations (\ref{N2}) and (\ref{N8}), we see that
a function $h$ belongs to $\Sigma [1]$ if and only if the function
$f$ defined by (\ref{N8}) belongs to the class $\mathcal{G}[1]$.
Moreover, $\left\vert \Im h(z)\right\vert \leq M$ for some
$M<\infty $ if and only if $f$ is of hyperbolic type (see
\cite{E-D-2006} and \cite{C-D-P}). In our setting this case
corresponds to the class $\Sigma _{A}^{\circ }[1]$ (i.e., $\alpha
=0$).

Now we begin to study the classes $\Sigma _{A}^{\alpha }[1]$,
respectively, $\Sigma ^{\alpha }[1]$, $\alpha >0$, of functions
$h\in \Sigma [1]$, for which
\begin{equation*}
\angle \lim_{z\rightarrow 1}(1-z)^{1+\alpha }h^{\prime }(z)=:\mu \text{,}
\end{equation*}
exists finitely.

\begin{remark} \label{remarkD}
Note that these classes arise naturally, if $h(\Delta )$ is
contained in a half-plane $\Pi =\{w\in \mathbb{C}:
w=h_{1}(z):=\dfrac{ibz}{1-z},$ $z\in \Delta \}$ and $\psi
=h_{1}^{-1}\circ h $ has a Dini smooth corner of opening $\pi
\alpha $, $0<\alpha \leq 1$, at the point $z=1$ (see \cite{PC-92},
p. 52). Indeed, in this case the limits
\begin{equation*}
\lim_{z\rightarrow 1}\frac{1-\psi (z)}{(1-z)^{\alpha }}\text{ \
and \ } \lim_{z\rightarrow 1}\frac{\psi ^{\prime
}(z)}{(1-z)^{\alpha -1}}\text{\ }
\end{equation*}
exist finitely and are different from zero. Moreover, it follows from Theorem 3.9
in \cite{PC-92}, that the second limit is exactly $\alpha $. Thus we have
from formula (\ref{n8''}) that the limit
\begin{eqnarray*}
\lim_{z\rightarrow 1}(1-z)^{1+\alpha }h^{\prime }(z)
&=&\lim_{z\rightarrow 1} \frac{(1-z)^{\alpha }}{1-\psi (z)}\cdot
\\ \cdot \lim_{z\rightarrow 1}\frac{(1-z)\psi ^{\prime
}(z)}{1-\psi (z)}\cdot ib &=&:\mu
\end{eqnarray*}
is also finite.
\end{remark}

Moreover, in this case formula (\ref{n9''}) implies that
\begin{equation*}
\nu :=\lim_{z\rightarrow 1}(1-z)^{\alpha
}h(z)=ib\lim_{z\rightarrow 1}\frac{ (1-z)^{\alpha }\psi
(z)}{1-\psi (z)}=\mu /\alpha
\end{equation*}
also exists finitely and is different from zero. So, $h(\Delta )$ has a corner
of opening $\pi \alpha $, $0<\alpha \leq 1$, at infinity (see \cite{PC-92},
p. 54).

We will see below (see Section 3) that if $h(\Delta )$ contains a
half-plane $h_{1}(\Delta )$ and $\varphi =h^{-1}\cdot h_{1}$ has a Dini
smooth corner of opening $\pi \gamma $, then $\frac{1}{2}\leq \gamma \leq 1$
and $h\in \Sigma _{A}^{\alpha }[1]$ with $\alpha =\frac{1}{\gamma }\in
\lbrack 1,2]$.

In general, we have the following simple assertion.

\begin{lemma} \label{lem3}
If $\Sigma _{A}^{\alpha }[1]\neq \varnothing $, then $\alpha
\in \lbrack 0,2]$.
\end{lemma}

\begin{proof} \label{proofN6}
Assume $\Sigma _{A}^{\alpha }[1]\neq \varnothing $ and let $h\in
\Sigma _{A}^{\alpha }[1]$. Once again, denote
\begin{equation*}
p(z):=\frac{1}{(1-z)^{2}h^{\prime }(z)}.
\end{equation*}%
Then we have $\Re p(z)\geq 0$ and
\begin{equation} \label{N12}
\angle \lim_{z\rightarrow 1}(1-z)^{\alpha -1}p(z)=\frac{1}{\mu }\neq
0,\infty .
\end{equation}%
But, as we have already seen,
\begin{equation*}
\angle \lim_{z\rightarrow 1}(1-z)p(z)\geq 0
\end{equation*}
is finite. This contradicts (\ref{N12}) if $\alpha >2$.
\end{proof}

To study the classes $\Sigma _{A}^{\alpha }[1]$ (respectively, $\Sigma
^{\alpha }[1]$) for $\alpha >0$, we need the following notion.

$\bullet $ We say that a generator $f$ belongs to the class
$\mathcal{G} _{A}^{\alpha }[1]$ ($\mathcal{G}^{\alpha }[1]$),
$\alpha >0$, if
\begin{equation} \label{*}
\angle \lim_{z\rightarrow 1}\frac{f(z)}{(1-z)^{1+\alpha }}=a\neq 0
\end{equation}
\begin{equation} \label{*'}
\left( \lim_{\substack{ z\rightarrow 1  \\ z\in \Delta
}}\frac{f(z)}{ (1-z)^{1+\alpha }}=a\neq 0\right) \text{.}
\end{equation}
It is clear that $h\in \Sigma _{A}^{\alpha }[1]$ (respectively, $\Sigma
_{A}^{\alpha }[1]$) if and only if the function $f$ defined by (\ref{N8}) belongs to the
class $\mathcal{G}_{A}^{\alpha }[1]$ (respectively, $\mathcal{G}^{\alpha
}[1] $) with $a=-\frac{1}{\mu }$. Of course, if $\alpha >0$, then $f$ is of
parabolic type.

The following example shows that the class $\Sigma ^{\alpha }[1]$
(respectively, $\mathcal{G}^{\alpha }[1]$) is a proper subset of
the class $\Sigma _{A}^{\alpha }[1]$ (respectively,
$\mathcal{G}_{A}^{\alpha }[1]$).

\begin{example} \label{exampleD}
Consider the function $f\in \Hol(\Delta , \mathbb{C})$ defined by
\begin{equation*}
f(z)=-(1-z)^{2}\left[ 1-\exp \left( -\frac{1+z}{1-z}\right) \right] ^{\beta
}(1-z)^{1-\beta }
\end{equation*}
$0<\beta \leq 1$.

It follows from the Berkson-Porta representation formula that $f\in \mathcal{G}[1]$.
Moreover, $f\in \mathcal{G}_{A}^{\alpha }[1]$, $\alpha =2-\beta $, since the
angular limit
\begin{equation*}
\angle \lim_{z\rightarrow 1}\frac{f(z)}{(1-z)^{1+\alpha }}=-1\text{.}
\end{equation*}%
At the same time the unrestricted limit
\begin{equation*}
\lim \frac{f(z)}{(1-z)^{1+\alpha }}
\end{equation*}%
does not exist, i.e., $f\notin \mathcal{G}^{\alpha }[1]$.
\end{example}

\begin{theorem} \label{teor3}
Let $h\in \Sigma _{A}^{\alpha }[1]$ with
\begin{equation*}
\angle \lim_{z\rightarrow 1}(1-z)^{1+\alpha }h^{\prime }(z)=:\mu .
\end{equation*}
The following two assertions hold:

(A) \ $\alpha \in \lbrack 0,2]$.

(B) \ Assume that for some $z\in \Delta $, the trajectory
$\{F_{t}(z)\}_{t\geq 0}$, where $F_{t}(z)=h^{-1}(h(z)+t)$,
 converges to $z=1$ nontangentially.

Then

(i) \ $\left\vert \arg \mu \right\vert \leq \dfrac{\pi }{2}\min
\{\alpha ,2-\alpha \}$. Moreover, if $\alpha \in (0,1]$, then this
inequality is sharp;

(ii) \ $\lim_{t\rightarrow \infty }t(1-F_{t}(z))^{\alpha
}=\dfrac{\mu }{ \alpha }$ for all $z\in \Delta $;

(iii) $\lim_{t\rightarrow \infty }\arg (1-F_{t}(z))=\dfrac{1}{\alpha }\arg
\mu $ for all $z\in \Delta $.
\end{theorem}

\bigskip To establish Theorem \ref{teor3}, we first need the following lemma.

\begin{lemma} \label{lem4}
Let $p\in \Hol(\Delta , \mathbb{C} )$ belong to the
Carath\'{e}odory class, i.e., $\Re p(z)\geq 0$ for all $z\in
\Delta $. Assume that for some $k\in \mathbb{R}$, the angular
limit
\begin{equation*}
\angle \lim_{z\rightarrow 1}(1-z)^{k}p(z)=\gamma
\end{equation*}
exists finitely and is different from zero. Then

(i) \ $k\in \lbrack -1,1];$

(ii) $\left\vert \arg \gamma \right\vert \leq \dfrac{\pi }{2}(1-\left\vert
k\right\vert )$.
\end{lemma}

\begin{proof} \label{proofN5}
Let $p\in \Hol(\Delta ,\Pi _{+})$ and assume that for some $k\in
\lbrack -1,1]$, the angular limit
\begin{equation*}
\angle \lim_{z\rightarrow 1}(1-z)^{k}p(z)=\gamma \neq 0.
\end{equation*}%
Define a function $\widehat{p}\in \Hol(\Pi _{+},\Pi _{+})$ by
\begin{equation*}
\widehat{p}(w)=p\left( \frac{w-1}{w+1}\right) \text{, }w\in \Pi _{+}.
\end{equation*}
Then we have
\begin{eqnarray} \label{NS}
\gamma &=&\angle \lim_{w\rightarrow \infty
}\widehat{p}(w)\frac{2^{k}}{ (w+1)^{k}}=   \\ &=&\angle
\lim_{w\rightarrow \infty }2^{k}\frac{\widehat{p}(w)}{w^{k}}.
\notag
\end{eqnarray}

Setting here $w=r\rightarrow \infty $, we obtain that
$\dfrac{\gamma }{2^{k}} \in \Pi ^{+}\backslash \{0\}$, so
\begin{equation*}
\left\vert \arg k\right\vert \leq \frac{\pi }{2}.
\end{equation*}
Now fix $\beta \in \left( -\dfrac{\pi }{2},\dfrac{\pi }{2}\right)
$ and set $w=re^{i\beta }\rightarrow \infty $ as $r\rightarrow
\infty $ in (\ref{NS}).

Then we have
\begin{equation*}
\frac{\gamma }{2^{k}}=\lim_{r\rightarrow \infty }\frac{\widehat{p}%
(re^{i\beta })}{(re^{i\beta })^{k}}=\lim_{r\rightarrow \infty }\frac{1}{r^{k}%
}\widehat{p}(re^{i\beta })e^{-i\beta k}.
\end{equation*}
Consequently,
\begin{equation*}
\arg \gamma =\lim_{r\rightarrow \infty }\left( \arg \widehat{p}(re^{i\beta
})-k\beta \right) .
\end{equation*}
At the same time, $\arg (\widehat{p}(re^{i\beta }))\in \left[
-\dfrac{\pi }{2 },\dfrac{\pi }{2}\right] $ for all $r>0$ and
$\beta \in \left( -\dfrac{\pi }{2},\dfrac{\pi }{2}\right) $.

Therefore
\begin{equation*}
\arg \gamma +k\beta =\lim_{r\rightarrow \infty }\left( \widehat{p}
(r(e^{i\beta }))\right) \in \left[ -\frac{\pi }{2},\frac{\pi
}{2}\right] .
\end{equation*}
Since $\beta \in \left( -\dfrac{\pi }{2},\dfrac{\pi }{2}\right) $
is arbitrary, we conclude that $\left\vert \arg \gamma \right\vert
\leq \dfrac{\pi }{2} (1-\left\vert k\right\vert )$.
\end{proof}

\begin{remark} \label{rem3}
The proof of Lemma 13 presented here is due to Santiago
D\'{\i}az-Madrigal \cite{DS}. This lemma as well as Theorem
\ref{teor3} were proved in \cite{E-S-Y} under the stronger
restriction that the unrestricted limits exist. So, Theorem
\ref{teor3} improves upon the results obtained there.
\end{remark}

\begin{proof12}
Since $f\left( =-\dfrac{1}{h^{\prime }} \right) \in
\mathcal{G}_{A}^{\alpha }[1]$ admits the representation (\ref{N6}),
assertion (A) of the theorem is a direct consequence of Lemma
\ref{lem4}. Assume now that there is a point $z\in \Delta $ such
that the trajectory $\{F_{t}(z)\}_{t\geq 0}$ converges to the
Denjoy-Wolff point $\tau =1$ nontangentially. For this point $z\in
\Delta $ we denote $u(t):=F_{t}(z)$, $t\geq 0$. Then it follows
from the Cauchy problem (\ref{N9}) that
\begin{equation*}
\begin{array}{c}
\dfrac{du}{(1-u)^{1+\alpha }}=-\dfrac{f(u)}{(1-u)^{1+\alpha }}dt
\\ u(0)=z,
\end{array}
\end{equation*}
which is equivalent to the integral equation
\begin{equation*}
\frac{1}{t}\cdot \frac{1}{(1-u(t))^{\alpha }}=-\frac{\alpha }{t}
\int\limits_{0}^{t}\frac{f(u(s))ds}{(1-u(s))^{1+\alpha
}}+\frac{1}{t}\cdot \frac{1}{(1-z)^{\alpha }}.
\end{equation*}
Now it follows from (\ref{*}) that
\begin{equation*}
\lim_{t\rightarrow \infty }t(1-u(t))^{\alpha }=-\frac{1}{\alpha a}=\frac{\mu
}{\alpha }
\end{equation*}
and
\begin{equation*}
\lim_{t\rightarrow \infty }\arg (1-u(t))=-\frac{1}{\alpha }\arg (-a)=\frac{1%
}{\alpha }\arg \mu \text{.}
\end{equation*}
This means that the trajectory $\{F_{t}(z)\}_{t\geq 0}$ has an
asymptote at its attractive point $\tau =1$. Next, it follows from
Theorem 2.9 in \cite{C-D-P} that for each $z\in \Delta $, the
trajectory $\{F_{t}(z)\}_{t\geq 0}$ has the same asymptote and all
the trajectories converge to $\tau =1$ nontangentially.

Repeating the above considerations for an arbitrary $z\in \Delta $, we
arrive at assertions (ii) and (iii) of the theorem.

Applying again Lemma \ref{lem4} and (\ref{N6}), we obtain
\begin{equation*}
\left\vert \arg (-a)\right\vert \leq \frac{\pi }{2}\min \{\alpha ,2-\alpha
\}.
\end{equation*}

But if $\alpha \in (0,1],$ then equality here is impossible, since otherwise
we see by (iii) that
\begin{equation*}
\lim_{t\rightarrow \infty }\arg (1-F_{t}(z))=\pm \frac{\pi }{2},
\end{equation*}
which contradicts the nontagential convergence. This proves assertion (i)
and we are done.
\end{proof12}

\begin{remark} \label{rem4}
We note that the conclusion of Theorem \ref{teor3} remains, of
course, valid if we replace condition (\ref{*}) by (\ref{*'}). Moreover, the
requirement of nontangential convergence in this case is not necessary. In
addition, formulas (i) and (iii) show that if $\alpha >1$ ($\alpha \leq 2$),
then all the trajectories must necessarily converge in a nontagential way.
\end{remark}

This proves the following assertion.

\begin{corollary} \label{colNov}
Let $f\in \mathcal{G}^{\alpha }[1]$. The following
assertions are equivalent:

(i) \ for some $z\in \Delta $, the trajectory $\{F_{t}(z)\}_{t\geq 0}$
converges to $z=1$ tangentially;

(ii) \ for all $z\in \Delta $, the trajectories $\{F_{t}(z)\}_{t\geq 0}$
converge to $z=1$ tangentially;

(iii) $\left\vert \arg (-a)\right\vert =\dfrac{\pi }{2}\alpha .$

Moreover, in this case $0<\alpha \leq 1$.
\end{corollary}

\begin{remark} \label{remarN}
Note that the nontangential convergence of the trajectory $
\{F_{t}(z)\}_{t\geq 0}$ to the Denjoy-Wolff point $z=1$ of the
semigroup $S$ means that
\begin{equation*}
\inf_{t\geq 0}\frac{1-\left\vert F_{t}(z)\right\vert }{\left\vert
1-F_{t}(z)\right\vert }>0.
\end{equation*}
We will see that this condition implies that there is no
horizontal half-plane which contains $h(\Delta )$.

In the opposite direction, the condition
\begin{equation*}
\lim_{t\geq 0}\frac{1-\left\vert F_{t}(z)\right\vert }{\left\vert
1-F_{t}(z)\right\vert }=0
\end{equation*}%
implies that the trajectory $\{F_{t}(z)\}_{t\geq 0}$ converges tangentially
to $z=1$ and conversely.

So, the question is whether the latter condition is sufficient to ensure
that $h(\Delta )$ lies in a horizontal half-plane.
\end{remark}

The following theorem shows that for $h\in \Sigma ^{\alpha }[1]$
the existence of a horizontal half-plane containing $h(\Delta )$
is, in fact, equivalent to the stronger condition that the ratio
$\dfrac{1-\left\vert F_{t}(z)\right\vert }{\left\vert
1-F_{t}(z)\right\vert }$ converges to zero faster than
$\dfrac{1}{t}$. Moreover, for $\alpha =1,$ we will see that the
limit
\begin{equation*}
\lim_{t\rightarrow \infty }\frac{t(1-\left\vert
F_{t}(z)\right\vert )}{ \left\vert 1-F_{t}(z)\right\vert }=L(z)
\end{equation*}
exists and is different from zero if and only if all the trajectories
converge to $z=1$ strongly tangentially, i.e., for each $z\in \Delta $ there
is a horocycle internally tangent to $\partial \Delta $ at the point $z=1$
which does not contain $\{F_{t}(z)\}_{t\geq s}$ for some $s\geq 0$.

We recall that by $\Sigma ^{\alpha }[1]$, $\alpha \in \lbrack
0,2]$, we denote the subclass of $\Sigma [1]$ which consists of
those functions $h$ for which the unrestricted $\lim_{\substack{
z\rightarrow 1  \\ z\in \Delta }}(1-z)^{1+\alpha }h^{\prime }(z)$
exists finitely.

\begin{theorem} \label{ter2}
Let $h\in \Sigma ^{\alpha }[1]$ with $\alpha \in (0,2]$,
\begin{equation*}
\mu :=\lim_{\substack{ z\rightarrow 1  \\ z\in \Delta }}(1-z)^{1+\alpha
}h^{\prime }(z)\neq \infty ,
\end{equation*}
and let $h^{-1}(h(z)+t)=:F_{t}(z)$.

Then $h(\Delta )$ lies in a horizontal half-plane if and only if
\begin{equation*}
\sup_{t\geq 0}\frac{t\left( 1-\left\vert F_{t}(z)\right\vert
\right) }{ \left\vert 1-F_{t}(z)\right\vert }=:M(z)<\infty
\text{.}
\end{equation*}

Moreover, in this case the following assertions hold:

(a) $\alpha $ must belong to the half-open interval $(0,1]$;

(b) \ $\left\vert \arg \mu \right\vert =\dfrac{\pi }{2}\alpha $;

(c) all the trajectories $\{F_{t}(z)\}$, $t\geq 0$, $z\in \Delta $, converge
tangentially to the point $z=1$.

(d) if $\Pi =\{w=\dfrac{ibz}{1-z}$: $z\in \Delta \}$ is a
half-plane containing $h(\Delta )$, then
\begin{equation*}
b\cdot \arg \mu >0\text{.}
\end{equation*}
\end{theorem}

\bigskip

In other words, $\Im h(z)$ is bounded from above if and only if
$\arg \mu >0$. Otherwise ($\arg \mu <0$), the set $\{\Im h(z):
z\in \Delta \}$ is bounded from below. Recall that since $\alpha
\neq 0$, $\Im h(z)$ cannot be bounded from both above and below.

\begin{corollary} \label{col1}
Let $h\in \Sigma ^{\alpha }[1]$ with $\alpha \in (0,2]$ and
\begin{equation*}
\mu :=\lim_{\substack{ z\rightarrow 1  \\ z\in \Delta }}(1-z)^{1+\alpha
}h^{\prime }(z)\neq \infty \text{.}
\end{equation*}
If either $\alpha >1$ or $\left\vert \arg \mu \right\vert \neq
\frac{\pi }{2} \alpha $, then there is no horizontal half-plane
which contains $h(\Delta )$.

In particular, if $\mu $ is real, then
\begin{equation*}
\sup_{z\in \Delta }\Im h(z)=\infty
\end{equation*}
and
\begin{equation*}
\inf_{z\in \Delta }\Im h(z)=-\infty \text{.}
\end{equation*}
\end{corollary}

\bigskip

\begin{example} \label{exam1}
Consider the function $h\in \Hol(\Delta , \mathbb{C})$ defined by
\begin{equation*}
h(z)=\frac{\mu}{K+1}\left[ \frac{1}{(1-z)^{K+1}}-1\right] , \quad z\in \Delta ,
\end{equation*}
with $K\neq -1$, $K\in \mathbb{R}$ and
$\mu\in\mathbb{C}\setminus\{0\}$.

Since $h(0)=0$ and
\begin{equation*}
h^{\prime }(z)=\frac{\mu}{(1-z)^{2+K}},
\end{equation*}
we have
\begin{equation*}
\Re \left[ (1-z)^{2}h^{\prime }(z)\right] =\Re \left[\mu\frac{1}{
(1-z)^{K}}\right] \geq 0
\end{equation*}
if and only if $-1<K\leq 1$ and $|\arg\mu|\leq \frac\pi2- \frac{\pi K}{2}$.

We also see that
\begin{equation*}
\angle \lim_{z\rightarrow 1}(1-z)^{1+\alpha }\cdot h^{\prime }(z)=\mu
\lim_{z\rightarrow 1}(1-z)^{1+\alpha }\cdot \frac{1}{(1-z)^{K+2}}
\end{equation*}
is finite and different from zero (and equal to $\mu $) if and only if $\alpha
=K+1$. Thus $\alpha $ must lie in $(0,2]$.

If $K>0$ or $|\arg\mu|<\frac\pi2- \frac{\pi K}{2}$, then it follows
from Theorem~ \ref{ter2} and Corollary~\ref{col1} that there is no
horizontal half-plane which contains $h(\Delta)$.
\end{example}

To prove Theorem \ref{ter2} we use again the relationship between
the class $ \Sigma [1]$ of functions convex in the positive
direction of the real axis and the class $\mathcal{G}[1]$ of
semigroup generators established in formulas (\ref{N2}) and
(\ref{N8}).

\begin{proof18}
For a domain $\Omega \subset \mathbb{C} $,  let $\delta
(w)=dist\left( w,\partial \Omega \right) $, $w\in \Omega $. Since
$h\in \Sigma ^{\alpha }[1]\subset \Sigma \lbrack 1]$, the fact
that $\Omega =h(\Delta )$ lies in a half-plane can be described as
follows: $\lim_{t\rightarrow \infty }\delta (h(z)+t)=:k(z)$ is
finite\ for each\ \ $ z\in \Delta $ (see \cite{C-D-P}).

Using again the Koebe Distortion Theorem, $\delta (h(\zeta ))\leq
(1-\left\vert \zeta \right\vert ^{2})\left\vert h^{\prime }(\zeta
)\right\vert \leq 4\delta (h(\zeta ))$, setting $\zeta =F_{t}(z)$, and using
the equality $h(F_{t}(z))=h(z)+t$, we get

\begin{equation*}
\delta (h(z)+t)\leq (1-\left\vert F_{t}(z)\right\vert ^{2})\left\vert
h^{\prime }(F_{t}(z))\right\vert \leq 4\delta (h(z)+t)
\end{equation*}
or
\begin{equation*}
\delta (h(z)+t)\leq \frac{t(1-\left\vert F_{t}(z)\right\vert
^{2})}{ \left\vert 1-F_{t}(z)\right\vert }\left\vert
1-F_{t}(z)\right\vert ^{1+\alpha }\left\vert h^{\prime
}(F_{t}(z))\right\vert \frac{1}{t\left\vert 1-F_{t}(z)\right\vert
^{\alpha }}\leq 4\delta (h(z)+t).
\end{equation*}

Since $h\in \Sigma ^{\alpha }[1]$, we see that for each $z\in \Delta $,

\begin{equation*}
\lim_{t\rightarrow \infty }\left\vert 1-F_{t}(z)\right\vert ^{1+\alpha
}\left\vert h^{\prime }(F_{t}(z))\right\vert =\left\vert \mu \right\vert .
\end{equation*}
In addition, we know that in this case the generator $f$ of the
semigroup $ h^{-1}(h(z)+t)=:F_{t}(z)$ belongs to
$\mathcal{G}^{\alpha }[1]$. Thus it follows from assertion (ii) of
Theorem \ref{teor3} and the remark following its proof that
$\lim_{t\rightarrow \infty }t\left\vert 1-F_{t}(z)\right\vert
^{\alpha }=\dfrac{1}{\alpha \left\vert \mu \right\vert }$ for all
$z\in \Delta $. Now recall that
\begin{equation*}
a:=\lim_{\substack{ z\rightarrow 1  \\ z\in \Delta }}\frac{f(z)}{
(1-z)^{1+\alpha }}=\frac{-1}{\mu }.
\end{equation*}
So the function $k(z)$ is finite if and only if for each $z\in \Delta $,
\begin{equation*}
\sup_{t\geq 0}\frac{t\left( 1-\left\vert F_{t}(z)\right\vert
\right) }{ \left\vert 1-F_{t}(z)\right\vert }=:M(z)<\infty
\text{.}
\end{equation*}
This condition, in its turn, implies that
\begin{equation*}
\lim_{t\rightarrow \infty }\frac{1-\left\vert F_{t}(z)\right\vert
}{ \left\vert 1-F_{t}(z)\right\vert }=0\text{,}
\end{equation*}
which means that all the trajectories $\{F_{t}(z)\}_{t\geq 0}$
converge tangentially to the point $z=1$ as $t\rightarrow \infty
$. Then, by Theorem \ref{teor3} (iii), $\left\vert \arg \mu
\right\vert =\dfrac{\pi }{2}\alpha $ . In addition, inequality (i)
in Theorem \ref{teor3} shows that if $\alpha
>1 $, then $\left\vert \arg \mu \right\vert =\dfrac{\pi }{2}\alpha \leq
\dfrac{\pi }{2}\left( 2-\alpha \right) $, which is impossible.

To prove (d), our last assertion, we consider the function $\psi
:\Delta \mapsto \Delta $ defined by $\psi (z)=h_{1}^{-1}(h(z))$,
where $h_{1}(z)= \dfrac{ibz}{1-z}$.

Then we have by Proposition 4 (see formulas (\ref{n4''}),
(\ref{n5''})).
\begin{eqnarray*}
\lim_{t\rightarrow \infty }\frac{1-\psi
(F_{t}(z))}{(1-F_{t}(z))^{\alpha }} &=&\lim_{t\rightarrow \infty
}\frac{t(1-G_{t}(\psi (z))}{t(1-F_{t}(z))^{ \alpha }}= \\
&=&\frac{\alpha }{\mu }\lim_{t\rightarrow \infty }t\left(
1-\frac{ib\psi (z)+t(1-\psi (z))}{ib+t(1-\psi (z))}\right)
=\frac{\alpha }{\mu }\cdot ib \text{.}
\end{eqnarray*}

If $\alpha <1$, then the function $\dfrac{1-\psi
(z)}{(1-z)^{\alpha }}$ does not admit negative real values. Then
it follows from the generalized Lindel\"{o}f Theorem (see Theorem
2.20 in \cite{G-L}, p. 42) that
\begin{equation} \label{d1}
\gamma =\angle \lim_{z\rightarrow 1}\frac{1-\psi
(z)}{(1-z)^{\alpha }}=\frac{ \alpha }{\mu }\cdot ib\text{.}
\end{equation}

If $\alpha =1$, then this limit also exists by the Julia-Wolff-Carath\'eodory
Theorem (see, for example, \cite{SD} and \cite{Shapiro}). Note that in this
case $\mu $ is pure imaginary; hence the limit $\gamma $ in (\ref{d1}) is a
real number.

In addition, since the function $1-\psi (z)$ is of positive real
part, we get by Lemma \ref{lem4} that
\begin{equation} \label{d2}
\left\vert \arg \gamma \right\vert \leq \frac{\pi }{2}(1-\alpha )\text{.}
\end{equation}

Then we have by (\ref{d1}) and (\ref{d2}),
\begin{eqnarray*}
\arg (ib) &=&\arg \gamma +\arg \mu > \\
&\geq &-\frac{\pi }{2}(1-\alpha )\geq -\frac{\pi }{2}
\end{eqnarray*}
if $\arg \mu >0$. That is, $b>0$. Similarly, if $\arg \mu <0$, then
$\arg (ib)\leq \dfrac{\pi }{2}(1-\alpha )<\dfrac{\pi }{2}$, that
is, $b<0$.

Noting that by assertion (ii), $\left\vert \arg \mu \right\vert
=\dfrac{\pi }{ 2}\alpha $ and $\alpha >0$, we conclude that $\mu $
cannot be a real number and Theorem 18 is proved.
\end{proof18}

Now, using Proposition \ref{Prop01}, we get the following assertion.

\begin{theorem} \label{teor.N3}
Let $S=\{F_{t}\}_{t\geq 0\text{ }}$ be a semigroup of parabolic
self-mappings of $\Delta $ generated by $f\in \mathcal{G}^{\alpha
}[1]$, i.e.,
\begin{equation*}
\lim_{z\rightarrow 1}\frac{f(z)}{(1-z)^{1+\alpha }}=a
\end{equation*}
exists finitely and is different from zero.

The following two assertions are equivalent:

(i) \ $\{F_{t}\}_{t\geq 0}$ outer conjugates with a group $\{G_{t}\}_{t\in
\mathbb{R} }$ of parabolic automorphisms of $\Delta $, i.e., there
is a conformal self-mapping $\psi :\Delta \mapsto \Delta $ such
that
\begin{equation*}
\psi (F_{t}(z))=G_{t}(\psi (z))\text{, \ \ }t\geq 0\text{, \ }z\in \Delta
\text{;}
\end{equation*}

(ii)
\begin{equation*}
\sup_{t\geq 0}\frac{t(1-\left\vert F_{t}(z)\right\vert )}{\left\vert
1-F_{t}(z)\right\vert }=:M(z)<\infty \text{, \ \ }z\in \Delta \text{.}
\end{equation*}

Moreover, in this case

(a) $\alpha \in (0,1]$;

(b) $\left\vert \arg (-a)\right\vert =\frac{\pi }{2}\alpha $;

(c) all the trajectories $\{F_{t}(z)\}_{t\geq 0}$, $z\in
\Delta $, $ t\geq 0$, converge tangentially to the point $z=1$.
\end{theorem}

As we have already mentioned, when $\alpha =1$ condition (ii)
is equivalent to the so-called strongly tangential convergence
(see Definition 22 below); hence it implies the existence of the
conjugating function.

Moreover, it turns out that strongly tangential convergence is possible
if and only if $\alpha =1$.

To see this, we denote by $d(z)$ the non-euclidean distance from a
point $ z\in \Delta $ to the boundary point $\zeta =1$:
\begin{equation*}
d(z):=\frac{\left\vert 1-z\right\vert ^{2}}{1-\left\vert
z\right\vert ^{2}} \text{.}
\end{equation*}
The sets $ \{ z\in\Delta : d(z)<k\}$, $k>0$, are horocycles
internally tangent to $\partial \Delta $ at $z=1$. It follows from
the Julia-Wolff-Carath\'eodory theorem for semigroups (see
\cite{E-S1} and \cite{SD}) that for each $z\in \Delta $, the
function $d(F_{t}(z))$ is decreasing, so the limit
\begin{equation*}
\lim_{t\rightarrow \infty }d(F_{t}(z))=\varepsilon (z)
\end{equation*}
exists and is nonnegative.

It can be shown that either $\varepsilon (z)>0$ for all $z\in
\Delta $ or $ \varepsilon (z)=0$ identically (see, for example,
\cite{SD-07}).

\begin{definition} \label{definN}
We say that a semigroup $S=\{F_{t}\}_{t\geq 0}$ is strongly
tangentially convergent if $\varepsilon (z)>0$.
\end{definition}

Geometrically, this means that for each $z\in \Delta $, there is a horocycle
$E$ internally tangent to $\partial \Delta $ at $z=1$ such that the
trajectory $\{F_{t}(z)\}_{t\geq 0}$ lies outside $E$.

Sometimes a strongly tangentially convergent semigroup is said to be of
finite shift (see \cite{C-M-P}).

\begin{theorem} \label{teor5}
Let $h\in \Sigma ^{\alpha }[1]$, $\alpha \in \lbrack 0,2]$, and let
$S=\{F_{t}\}_{t\geq 0}$ be defined by $F_{t}(z)=h^{-1}(h(z)+t)$, $
z\in \Delta $, $t\geq 0$. The following two assertions hold.

(A) Assume that for some $z\in \Delta $, the trajectory $\{F_{t}(z)\}_{t\geq
0}$ converges to $z=1$ strongly tangentially. Then

\ \ \ \ (i) $\ \alpha =1$;

\ \ \ (ii) $h(\Delta )$ is contained in a half-plane.

(B) Conversely, assume that conditions (i) and (ii) are fulfilled.
Then all the trajectories $\{F_{t}(z)\}_{t\geq 0}$, $z\in \Delta
$, converge to $z=1$ strongly tangentially.

Moreover,
\begin{equation*}
\mu =\lim_{z\rightarrow 1}(1-z)^{2}h^{\prime }(z)
\end{equation*}
is purely imaginary, and either

$\Im h(z)\geq \dfrac{-\Im \mu }{2\varepsilon (0)},$ $z\in \Delta $
, if $\Im \mu >0$

or

$\Im h(z)\leq \dfrac{-\Im \mu }{2\varepsilon (0)}$, if $\Im \mu
<0$.
\end{theorem}

\begin{proof} \label{profN}
As we have already mentioned, if the trajectory
$\{F_{t}(z)\}_{t\geq 0}$ converges to $z=1$ (strongly)
tangentially, then $ 0<\alpha \leq 1$. We claim that $\alpha $ cannot
be less than $1$. Indeed, consider
\begin{equation*}
d(F_{t}(z)):=\frac{\left\vert 1-F_{t}(z)\right\vert ^{2}}{1-\left\vert
F_{t}(z)\right\vert ^{2}}\text{.}
\end{equation*}
Since $h\in \Sigma ^{\alpha }[1]$, we have that for $t$ large enough,
\begin{equation*}
\left\vert 1-F_{t}(z)\right\vert ^{1+\alpha }\left\vert h^{\prime
}(F_{t}(z)\right\vert \leq \left\vert \mu \right\vert +1\text{.}
\end{equation*}

Then again by the Koebe Distortion Theorem we get
\begin{eqnarray*}
dF_{t}(z) &=&\frac{\left\vert 1-F_{t}(z)\right\vert ^{1+\alpha }}{
(1-\left\vert F_{t}(z)\right\vert ^{2})}\frac{\left\vert h^{\prime
}(F_{t}(z))\right\vert }{\left\vert h^{\prime
}(F_{t}(z))\right\vert }\cdot \left\vert 1-F_{t}(z)\right\vert
^{1-\alpha }\leq \\ &\leq &\frac{\left\vert \mu \right\vert
+1}{\delta (h(z)+t)}\left\vert 1-F_{t}(z)\right\vert ^{1-\alpha
}\text{.}
\end{eqnarray*}
Thus, if at least one of conditions (i) or (ii) does not hold,
i.e., $\alpha <1$ or $\lim_{t\rightarrow \infty }\delta
(h(z)+t)=\infty $, then $ d(F_{t}(z))$ must converge to zero, and
this is a contradiction. Assertion (A) is proved.

Now we prove the converse assertion (B). To this end, assume that (i) and (ii) hold.

Fix any $w\in h(\Delta )$ and denote by $\delta (w)$ the euclidean distance
from $w$ to $\partial h(\Delta )$.

Then condition (ii) of the theorem can be rewritten as follows:
\begin{equation} \label{N2'}
\lim_{t\rightarrow \infty }\delta (h(z)+t)=k(z),\text{ \ \ \ \ }z\in \Delta ,
\end{equation}
is finite.

Note that condition (\ref{N2'}) and Theorem 2.9 (b) in \cite{C-D-P} imply
that all trajectories converge to the point $\tau =1$ tangentially. To prove
strongly tangential convergence, we use again the Koebe Distortion Theorem
which asserts that for all $w\in h(\Delta )$,
\begin{equation}  \label{KT}
\delta (w)\leq (1-\left\vert \zeta \right\vert ^{2})\left\vert
h^{\prime }(\zeta )\right\vert \leq 4\delta (w),
\end{equation}
where $w=h(\zeta )$, $\zeta \in \Delta $.

Setting $\zeta =F_{t}(z)$ in (\ref{KT}) and using the equality $
h(F_{t}(z))=h(z)+t$, we have
\begin{equation*}
\delta (h(z)+t)\leq (1-\left\vert F_{t}(z)\right\vert ^{2})h^{\prime
}(F_{t}(z))\leq 4\delta (h(z)+t)
\end{equation*}
or
\begin{equation*}
\delta (h(z)+t)\leq \frac{\left\vert 1-F_{t}(z)\right\vert ^{2}\left\vert
h^{\prime }(F_{t}(z))\right\vert }{d(F_{t}(z))}\leq 4\delta (h(z)+t),
\end{equation*}
where $d(F_{t}(z))=\dfrac{\left\vert 1-F_{t}(z)\right\vert ^{2}}{
1-\left\vert F_{t}(z)\right\vert ^{2}}$ converges to $\varepsilon
(z)\geq 0$ as $t\rightarrow \infty $.

Since
\begin{equation*}
\lim_{t\rightarrow \infty }\left\vert 1-F_{t}(z)\right\vert ^{2}\left\vert
h^{\prime }(F_{t}(z))\right\vert =\left\vert \mu \right\vert ,
\end{equation*}
we get the inequality
\begin{equation*}
k(z)\leq \frac{\left\vert \mu \right\vert }{\varepsilon (z)}\leq 4k(z),
\end{equation*}
which shows that $k(z)$ is finite if and only if $\varepsilon (z)>0$. This
proves the implication (B)$\Longrightarrow $(A).

To prove our last assertion, we use Theorem 3 from \cite{SD-07}
which asserts that if $f\in \mathcal{G}^{1}[1]$ with $f^{\prime
\prime }(1)\neq 0$, then the trajectories $\{F_{t}(z)\}_{t\geq
0}$, $z\in \Delta $, converge to $\tau =1$ strongly tangentially
if and only if the net
\begin{equation} \label{N3'}
g_{t}(z)=d_{t}(0)\left[ \frac{1+F_{t}(z)}{1-F_{t}(z)}-\frac{2i\Im
\overline{F_{t}(0)}}{\left\vert 1-F_{t}(0)\right\vert ^{2}}\right]
\end{equation}
converges locally uniformly as $t\rightarrow\infty$ to a holomorphic function
\begin{equation*}
g:\Delta \mapsto \Pi _{+}=\left\{ w\in \mathbb{C} :\Re w\geq
0\right\} \text{,}
\end{equation*}
which solves Abel's equation
\begin{equation}  \label{N4'}
g(F_{t}(z))=g(z)+ibt\text{,}
\end{equation}
where $b=-\varepsilon (0)\Im f^{\prime \prime }(1)$.

It was also shown there that the condition $\varepsilon (0)>0$
implies that $ a=\frac{1}{2}f^{\prime \prime }(1)$ (hence $\mu
=-\dfrac{1}{a}$) is purely imaginary.

Differentiating now (\ref{N4'}) at $t=0^{+}$, we get
\begin{equation*}
g^{\prime }(z)=-\frac{ib}{f(z)}.
\end{equation*}
In addition, (\ref{N3'}) implies that $\Re g(z)\geq 0$ and
$g(0)=1$. But it follows now from (\ref{N8}) that
\begin{equation*}
g(z)=-ibh(z)+1
\end{equation*}
or
\begin{eqnarray*}
h(z) &=&\frac{i}{-\varepsilon (0)\Im f^{\prime \prime }(1)}\left[ g(z)-1%
\right] = \\ &=&\frac{i\Im \mu }{2\varepsilon (0)}\left[
g(z)-1\right] .
\end{eqnarray*}

So, $\Im h(z)\geq -\Im \mu /2\varepsilon (0)$ if $\Im \mu >0$ and
$\Im h(z)\leq -\Im \mu /2\varepsilon (0)$ if $\Im \mu <0$ .
Theorem 23 is proved.
\end{proof}

Alternatively, the implication (B)$\Longrightarrow $(A) can also
be proved by using the following Proposition and Theorem
\ref{ter2}.

\begin{proposition} \label{propoz}
Let $S=\{F_{t}\}_{t\geq 0\text{ }}$be a semigroup of parabolic
type generated by $f\in \mathcal{G}^{1}[1]$, i.e.,
\begin{equation*}
\lim_{z\rightarrow 1}\frac{f(z)}{(1-z)^{2}}=:\frac{1}{2}f^{\prime \prime
}(1)=a
\end{equation*}
exists finitely and is different from zero.

Then $\{F_{t}\}_{t\geq 0\text{ }}$ converges to $z=1$ strongly
tangentially if and only if
\begin{equation*}
\lim_{t\rightarrow \infty }\frac{t(1-\left\vert
F_{t}(z)\right\vert )}{ \left\vert 1-F_{t}(z)\right\vert }=L(z)
\end{equation*}
is finite.
\end{proposition}

\begin{proof} \label{proofN1}
We know already that
\begin{equation*}
\lim_{t\rightarrow \infty }t(1-F_{t}(z))=-\frac{1}{a}\text{.}
\end{equation*}
Therefore we have
\begin{eqnarray*}
\varepsilon (z) &=&\lim_{t\rightarrow \infty }\frac{\left\vert
1-F_{t}(z)\right\vert ^{2}}{1-\left\vert F_{t}(z)\right\vert
^{2}}= \\ &=&\frac{1}{2}\lim_{t\rightarrow \infty
}\frac{t\left\vert 1-F_{t}(z)\right\vert \cdot \left\vert
1-F_{t}(z)\right\vert }{ t(1-\left\vert F_{t}(z)\right\vert )}= \\
&=&\frac{1}{2\left\vert a\right\vert }\lim_{t\rightarrow \infty
}\frac{ \left\vert 1-F_{t}(z)\right\vert }{t(1-\left\vert
F_{t}(z)\right\vert )}=- \frac{1}{2\left\vert a\right\vert
L(z)}\text{.}
\end{eqnarray*}
Thus we see that $\varepsilon (z)>0$ if and only if $L(z)<\infty $.
\end{proof}

Using this theorem and Theorem \ref{ter2}, one can easily construct an
example of a semigroup which converges tangentially, but not strongly tangentially, to its Denjoy-Wolff
point.

\begin{example} \label{exam.N}
Consider a function $h:\Delta \mapsto \mathbb{C}$ conformally
mapping the open unit disk $\Delta $ onto the quadrant
$\{(x,y)\in\mathbb{R}^{2}: -1\leq x<\infty $, $-1\leq y<\infty \}$
given by the formula
\begin{equation*}
h(z)=e^{i\frac{\pi }{4}}\left[ \sqrt{\frac{1+z}{1-z}}-1\right] \text{.}
\end{equation*}
It is clear that $h\in \Sigma ^{\alpha }[1]$, with $\alpha =\frac{1}{2}$ and
\begin{equation*}
\mu =\lim (1-z)^{\alpha +1}h^{\prime }(z)=e^{i\frac{\pi }{4}}\text{.}
\end{equation*}
The corresponding generator $f\left( =-\frac{1}{h^{\prime }}\right) $ is
defined by $f(z)=-(1-z)^{2}\sqrt{\frac{1+z}{1-z}}e^{-\frac{\pi }{4}i}$.

Since $h(\Delta )$ lies in a horizontal half-plane, all the
trajectories $ F_{t}(z)$ must converge to $z=1$ tangentially:
\begin{equation*}
\lim_{t\rightarrow \infty }\arg (1-F_{t}(z))=\pm \frac{\arg \mu
}{\alpha } =\pm \frac{\pi }{2}\text{.}
\end{equation*}
At the same time no trajectory can converge to $z=1$ strongly
tangentially because $\alpha \neq 1$.
\end{example}

Combining this theorem with a result in \cite{SD-08}, we get the following
rigidity assertion.

\begin{corollary} \label{coll.N}
Let $f\in \mathcal{G}^{\alpha }[1]$ for some
$\alpha \in \lbrack 0,2]$ admit the representation
\begin{equation*}
f(z)=a(1-z)^{\alpha +1}+b(1-z)^{3}+\gamma (z),
\end{equation*}
where
\begin{equation*}
\angle \lim_{z\rightarrow 1}\frac{\gamma (z)}{(1-z)^{3}}=0\text{.}
\end{equation*}

Assume that the semigroup $S=\{F_{t}\}_{t\geq 0\text{ }}$ converges to $z=1$
strongly tangentially. Then

(i) \ $\alpha =1$; hence $\Re a=0$;

(ii) $\Im b=0$ and $\Re b\leq 0$;

(iii) $\Re b=0$ if and only if $\gamma \equiv 0$, i.e., $f$ is the
generator of a group of parabolic automorphisms of $\Delta $.
\end{corollary}

Under some additional conditions of regularity one can establish a more
general criterion for $h(\Delta )$ to lie in a half-plane (not necessarily
horizontal).

Let us assume that $h(\Delta )$ is a Jordan domain (the region
interior to a Jordan curve). Then the Carath\'{e}odory Extension
Theorem (see, for example, \cite{PC-92}) asserts that $h$ extends
to a homeomorphism of $\overline{\Delta }$ and formula (2) then
yields an extension of $F_{t}$ to a homeomorphism of
$\overline{\Delta }$ for each $t>0$.

We also assume that $h(z)+1\in h(\Delta )$ for all $z\in
\overline{\Delta }$ , which implies that $\overline{F_{1}(\Delta
)}\subset \Delta $.

Finally, we suppose that the generator $f:\Delta \mapsto
\mathbb{C} $ of the semigroup $S=\{F_{t}\}_{t\geq 0}$
$(f(z)=-1/h^{\prime }(z))$ admits the representation
\begin{equation*}
f(z)=a(z-1)^{2}+b(z-1)^{3}+\gamma (z),
\end{equation*}
where
\begin{equation*}
\lim_{z\rightarrow 1}\frac{\gamma (z)}{(z-1)^{3+\varepsilon }}=0
\end{equation*}
for some $\varepsilon >0$.

In this case one can show (see, for example, \cite{E-L-R-S}) that
the single mapping $F_{1}$ can be represented in the form
\begin{equation*}
F_{1}(z)=C^{-1}(w+a_{1}+\frac{b_{1}}{w+1}+\Gamma (w+1)),
\end{equation*}
where $w=C(z):=\dfrac{1+z}{1-z}$, $a_{1}=F^{\prime \prime }(1)$, $
b_{1}=(F^{\prime \prime }(1))^{2}-\dfrac{2}{3}F^{\prime \prime
\prime }(1)$ (the Schwarz derivative), and
\begin{equation*}
\left\vert \Gamma (w+1)\right\vert \leq \frac{C}{\left\vert w+1\right\vert
^{1+\varepsilon }}.
\end{equation*}

It was shown in \cite{BP-SJ} (p. 62) that in this case there is a
(unique) function $g:\Delta \mapsto \mathbb{C} $ with $g(0)=0$
which solves Abel's equation
\begin{equation*}
g(F_{1}(z))=g(z)+a_{1}\text{.}
\end{equation*}
Moreover, $\Re g(z)$ is bounded from below if and only if
\begin{equation}  \label{NC}
\Re \overline{a_{1}}\cdot b_{1}\leq 0\text{.}
\end{equation}
Using now Theorem 1 in \cite{E-L-R-S}, we obtain
\begin{equation*}
F^{\prime \prime }(1)=-f^{\prime \prime }(1)
\end{equation*}
and
\begin{equation*}
F^{\prime \prime \prime }(1)\equiv \frac{3}{2}(f^{\prime \prime
}(1))^{2}-f^{\prime \prime \prime }(1).
\end{equation*}
So, we see that
\begin{equation*}
a_{1}=F^{\prime \prime }(1)=-f^{\prime \prime }(1)=-2a
\end{equation*}
and
\begin{eqnarray*}
b_{1} &=&(F^{\prime \prime }(1))^{2}-\frac{2}{3}F^{\prime \prime \prime }(1)=
\\
&=&(F^{\prime \prime }(1))^{2}-\frac{2}{3}\left[ \frac{3}{2}(f^{\prime
\prime }(1))^{2}-f^{\prime \prime \prime }(1)\right] = \\
&=&\frac{2}{3}f^{\prime \prime \prime }(1)=\frac{2}{3}6b=4b.
\end{eqnarray*}

Thus condition (\ref{NC}) can be replaced by $\Re
\overline{a}b\leq 0$. Since $h(z)=\frac{1}{a}g(z)$, we get the
following assertion.

\begin{theorem} \label{teor6}
Let $h\in \Sigma ^{1}[1]$ be such that

(a) $h(\Delta )$ is a Jordan domain;

(b) $h(z)+1$ belongs to $h(\Delta )$ for all $z\in \overline{\Delta }$;

(c) $f(z)=-\frac{1}{h^{\prime }(z)}$ admits the representation
\begin{equation*}
f(z)=a(z-1)^{2}+b(z-1)^{3}+\gamma (z),
\end{equation*}
where
\begin{equation*}
\lim_{z\rightarrow 1}\frac{\gamma (z)}{(z-1)^{3+\varepsilon }}=0
\end{equation*}
for some $\varepsilon >0$.

Then $h( \Delta )$ lies in a half-plane if and only if
\begin{equation} \label{NC2}
\Re \overline{a}b\leq 0.
\end{equation}
\end{theorem}

\bigskip

\begin{remark} \label{rem2}
It follows immediately from the Berkson-Porta representation formula that $\Re
a\leq 0$. Moreover, under our assumption $\Re a=0$ if and only if
the trajectories $\{F_{t}(z)\}_{t\geq 0}$ converge to $z=1$
strongly tangentially.

In addition, since
\begin{equation*}
b=\lim_{z\rightarrow 1}\frac{q(z)}{1-z},
\end{equation*}
where
\begin{equation*}
q(z)=p(z)+a=-\frac{f(z)}{(1-z)^{2}}+a
\end{equation*}
is a function with a non-negative real part, we see that $\Re a=0$
implies that $b$ is a real number, so condition (\ref{NC2}) holds
automatically.
\end{remark}

Finally, we note that it was proved in \cite{SD-07} that if $\Re
a=b=0$ , then, in fact, $S=\{F_{t}\}_{t\geq 0}$ consists of
parabolic automorphisms of $\Delta $. In this case $h\in \Sigma
^{1}[1]$ must be of the form
\begin{equation*}
h(z)=\frac{ikz}{1-z}
\end{equation*}
for some real $k\neq 0$.

\section{Covering theorems (inner conjugation)}

\subsection{Backward flow invariant domains}

\ \ \ Let $f\in \mathcal{G}[1]$ be the generator of a semigroup $
S=\{F_{t}\}_{t\geq 0\text{ }}$ having the Denjoy-Wolff point $\tau
=1$, i.e.,
\begin{equation*}
f(1)\left( :=\angle \lim_{z\rightarrow 1}f(z)\right) =0
\end{equation*}
and
\begin{equation*}
f^{\prime }(1)\left( =\angle \lim_{z\rightarrow 1}\frac{f(z)}{z-1}\right)
\geq 0\text{.}
\end{equation*}

\bigskip

\begin{definition}
We say that a simply connected domain $\Omega \subset \Delta $
is a backward flow invariant domain (BFID) for $S$ if $S$ can be
extended to a group on $\Omega $.
\end{definition}

Since $\Omega $ (if it exists) is simply connected, one can find a
Riemann conformal mapping $\varphi :\Delta \mapsto \Delta $ such
that $\varphi (\Delta )=\Omega $. It is clear that the family
$\{G_{t}\}_{t\in \mathbb{R}}$ defined by
\begin{equation} \label{N1''}
G_{t}(z)=\varphi ^{-1}(F_{t}(\varphi (z)))\text{, \ \ }t\in
\mathbb{R} \text{,}
\end{equation}
forms a group of automorphisms of $\Delta $.

Thus, the existence of a BFID for $S$ is equivalent to the existence of
an inner conjugator $\varphi :\Delta \mapsto \Delta $,
\begin{equation} \label{N2''}
F_{t}(\varphi (z))=\varphi (G_{t}(z))\text{, \ \ \ }t\geq 0,
\end{equation}
where $\{G_{t}\}_{t\in \mathbb{R} }$ is a group on $\Delta $.

Furthermore, since $\tau =1$ is the Denjoy-Wolff point of $S$ it
follows from (\ref{N2''}), that
\begin{equation}  \label{N3''}
\lim_{t\rightarrow \infty }\varphi (G_{t}(z))=1\text{.}
\end{equation}
Hence, $\tau =1\in \overline{\Omega }$.

It is also clear clear that the group $\{G_{t}\}$ is not elliptic.
Indeed, otherwise there would be a point $\eta \in \Delta $ such that
$G_{t}(\eta )=\eta $ for all $t\in \mathbb{R} $. Then (\ref{N2})
would imply that $\varphi (\eta )\in \Delta $ is a fixed point of
$S:\varphi (\eta )=\varphi (G_{t}(\eta ))=F_{t}(\varphi (\eta ))$.
This would contradict our assumption that $S$ has a boundary
Denjoy-Wolff point.

Thus the group $\{G_{t}\}_{t\in \mathbb{R} }$ is either of
parabolic or hyperbolic type.

\begin{definition}
We say that $\Omega $ is of $p$-type if $\{G_{t}\}_{t\in
\mathbb{R} }$ is parabolic and that it is of $h$-type if this group is
hyperbolic.
\end{definition}

\begin{remark}
We will see below that if $\Omega $ is of $p$-type, then $
S=\{F_{t}\}_{t\geq 0\text{ }}$ must also be of parabolic type,
i.e., $ f^{\prime }(1)=0$. At the same time, if $\Omega $ is of
$h$-type then the semigroup $S$ might be of hyperbolic as well as
parabolic type.
\end{remark}

\begin{example}
Consider the holomorphic function $f_{1}$ defined by
\begin{equation*}
f_{1}(z)=-(1-z)^{2}\cdot \frac{2+\sqrt{\frac{1+z}{1-z}}}{1+\sqrt{\frac{1+z}{
1-z}}}\cdot \frac{1+z}{1-z}.
\end{equation*}
One can check by using the Berkson--Porta representation formula that $f_{1}\in
\mathcal{G} [1]$. Moreover, $f_{1}^{\prime }(1)=2$, i.e., the
semigroup $S_{1}$ generated by $f_{1}$ is of hyperbolic type. In
addition, $f_{1}(-1)=0$ and $ f_{1}^{\prime }(-1)=-4$. So, by
Theorem 1 in \cite{E-S-Z1}, $S_{1}$ has a BFID. Using Theorem~3 in
\cite{E-S-Z1}, one can find a univalent function $ \varphi $ which
maps $\Delta $ onto the maximal BFID $\Omega $ conformally. Namely,
\begin{equation*}
\varphi (z)=\frac{(w-1)^{2}-1}{(w-1)^{2}+1}\,,\quad
\mbox{where}\quad w= \sqrt{1+\sqrt{\frac{1+z}{1-z}}}\,.
\end{equation*}
In Figure~(a) we present the vector field of the flow $S_{1}$, as well as the
location of the h-type maximal BFID for $S_{1}$ in the unit disk $\Delta $.
\begin{figure}\centering \subfigure[h-type maximal BFID]{
    \includegraphics[angle=270,width=5.4cm,totalheight=5.4cm]{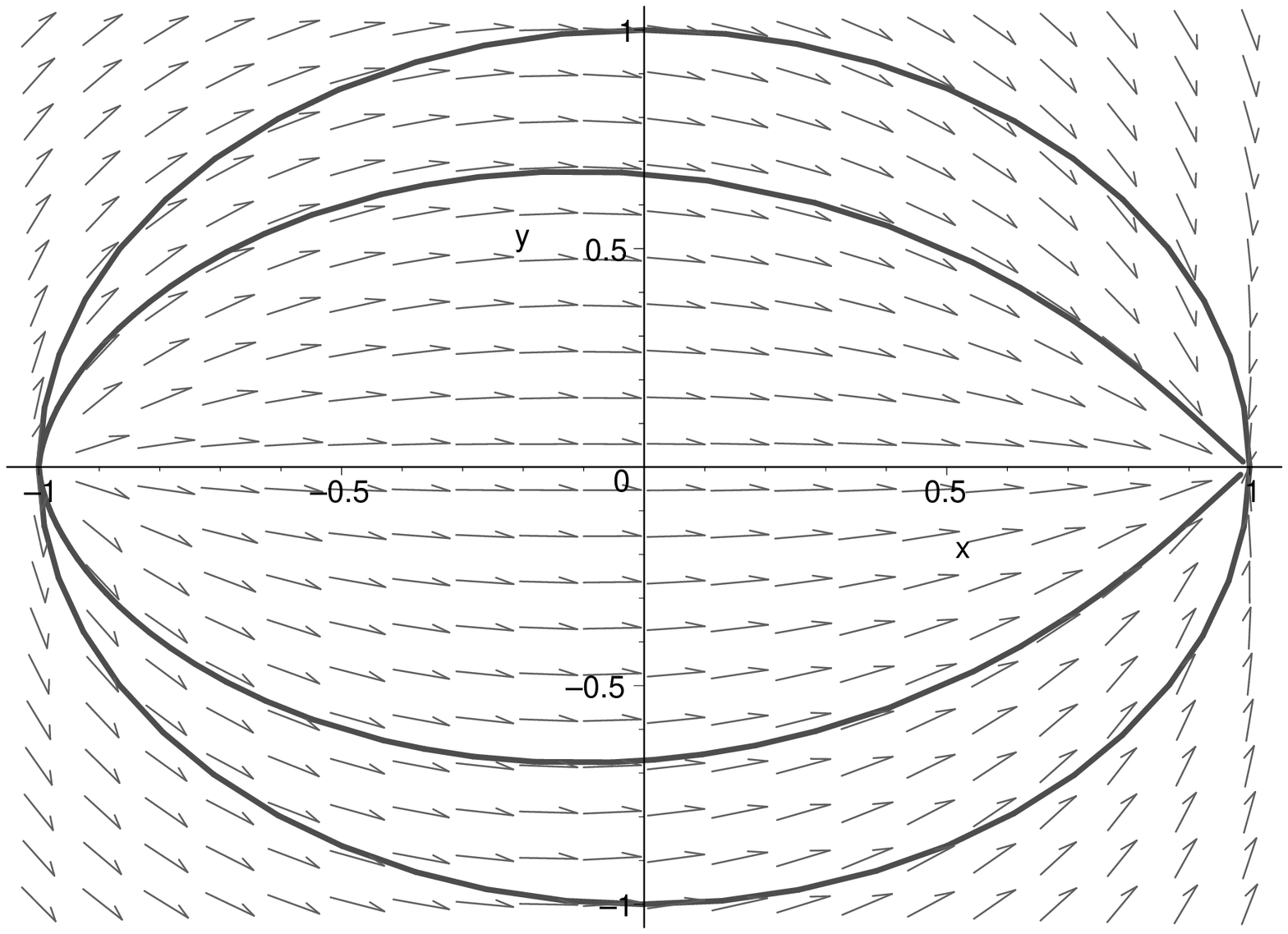}}
    \hspace{7mm}
\subfigure[BFIDs of h-type and p-type]{
    \includegraphics[angle=270,width=5.4cm,totalheight=5.4cm]{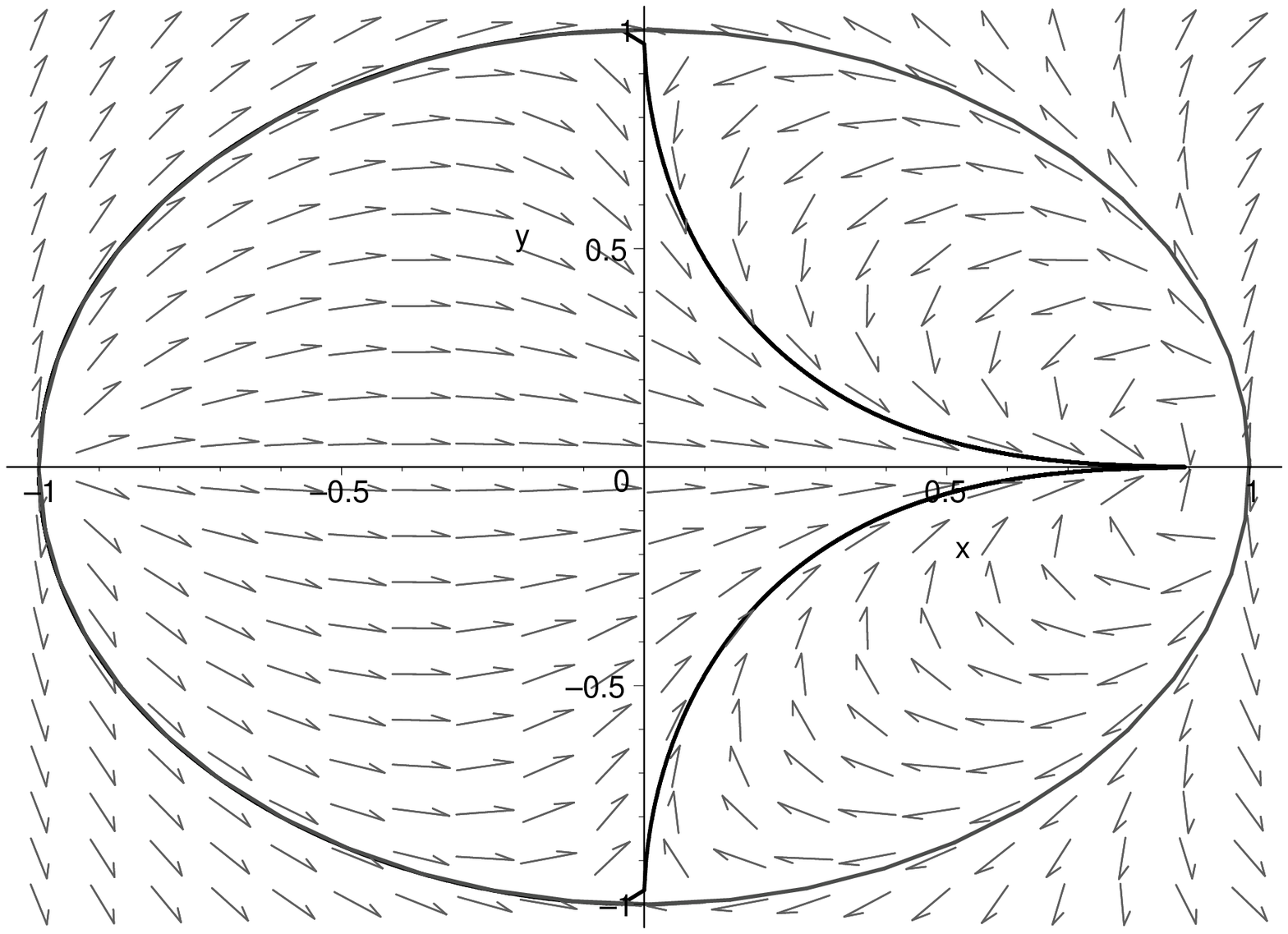}}
\end{figure}
\end{example}

\begin{example}
Consider now another semigroup generator:
\begin{equation*}
f_{2}(z)=-(1-z)^{2}\frac{1-z^{2}}{1+z^{2}}.
\end{equation*}
Since $f_{2}^{\prime }(1)=0$, we conclude that the semigroup
$S_{2}=\left\{ F_{t}\right\} _{t\geq 0}$ generated by $f_{2}$ is
of parabolic type. Using Theorems 1 and 3 in \cite{A-E-R-S}, and
Theorem~\ref{teoremaN} and Remark~\ref {remarkN} below, one can
show that $S_{2}$ has BFIDs of both types: two domains of p-type
and a domain of h-type. More precisely, in this case the closure
of the union of these three BFIDs is the whole disk $\Delta $. In
other words, there are only two trajectories $\left\{
F_{t}(z_{i})\right\} _{t\geq 0},\ i=1,2,$ that cannot be extended
to all the values of $t\in \mathbb{R}$. In Figure~(b) we show the
maximal BFIDs of both types.
\end{example}

\begin{theorem} \label{teoremaN}
Let $f\in \mathcal{G}[1]$ and let $S=\{F_{t}\}_{t\geq 0 \text{ }}$
be the semigroup generated by $f$. The following assertions are equivalent:

(i) there is a point $z\in \Delta $ such that $F_{t}(z)\in \Delta
$ for all $ t\in \mathbb{R} $ and the point $\zeta
=\lim_{t\rightarrow -\infty }F_{t}(z)$ is a boundary regular null
point of $f$, i.e., $f^{\prime }\left( \zeta \right) $ $=\angle
\lim_{z\rightarrow \zeta }f(z)$ is finite;

(ii) the image $h(\Delta ),$ where $h$ $\in \Sigma [1]$ is defined
by
\begin{equation} \label{D1}
h\left( z\right) =-\int\limits_{0}^{z}\frac{dz}{f\left( z\right)
},
\end{equation}
contains a horizontal strip (or a horizontal half-plane);

(iii) there is a simply connected open subset $\Omega \subset
\Delta $ such that $F_{t}|_{\Omega }$ is a group of automorphisms
of $\Omega $;

(iv) there are real numbers $a$ and $b$ such that the differential
equation
\begin{equation} \label{a*}
\left[ a(z^{2}-1)+ib(1-z)^{2}\right] \varphi ^{\prime }(z)=f(\varphi (z))
\end{equation}
has a nonconstant solution $\varphi \in \Hol(\Delta )$.

Moreover, in this case $\varphi $ is univalent and $\Omega =\varphi (\Delta
) $ is a BFID for $S$.
\end{theorem}

\begin{remark} \label{remarkN}
Actually, one can distinguish between the following two situations

(A) The following assertions are equivalent:

\ \ \ (a) there is a boundary regular point $\zeta \in \partial \Delta $ of the
\ generator $f$ which is different from the Denjoy-Wolff point $\tau =1$;

\ \ \ (b) \ $h(\Delta )$ contains a horizontal strip;

\ \ \ (c) there are numbers $a\neq 0$ and $b\in \mathbb{R} $ such
that equation (\ref{a*}) has a solution $\varphi $ for which
$\varphi (\Delta )=\Omega $ is a BFID of $h$-type, i.e.,
$S=\{F_{t}\}$ is a conjugate of a group of hyperbolic
automorphisms and $\partial \Omega \cap \partial \Delta \ni \zeta
\neq 1$.

(B) The following assertions are equivalent:

\ \ \ (a') there is a point $z\in \Delta $ such that
$\{F_{t}(z)\}$ is well defined for all $t\in \mathbb{R} $ and
\begin{equation*}
\lim_{t\rightarrow -\infty }F_{t}(z)=1\text{;}
\end{equation*}

\ \ \ (b') \ $h(\Delta )$ contains a half-plane;

\ \ \ (c') equation (\ref{a*}) with $a=0$ and $b\neq 0$ has a
nonconstant solution $\varphi $ such that $\varphi (\Delta
)=\Omega $ is a BFID of $p$ -type, i.e., $S=\{F_{t}\}$ is a
conjugate of a group of parabolic automorphisms and $\partial
\Omega \cap \partial \Delta \ni 1$.
\end{remark}

\begin{proof39}
The equivalence of assertions (i) and (ii) is shown in
\cite{C-D-P} and \cite{E-S-Z1}.

Furthermore, if $k$ is a Riemann mapping of $\Delta $ onto a
horizontal strip (or a half-plane) contained in $h(\Delta ),$ then
the function $ \varphi =h^{-1}\circ k$ is a univalent self-mapping
of $\Delta $ and $\Omega =\varphi (\Delta )$ is a BFID for $S$.
Indeed, the family $\{G_{t}\}_{t\in \mathbb{R} }$ defined by
$G_{t}\left( z\right) =k^{-1}(k\left( z\right) +t),$ $t\in $ $
\mathbb{R} ,$ forms a one-parameter group of hyperbolic ( or
parabolic) automorphisms of $\Delta$. On the other hand,

\begin{equation*}
\varphi \left( G_{t}\left( z\right) \right) =h^{-1}(k\left( z\right)
+t)=h^{-1}(h\left( \varphi \left( z\right) \right) +t)=F_{t}\left( \varphi
\left( z\right) \right) .
\end{equation*}
That is, $\varphi $ is, in fact, an inner conjugating function for
$ \{F_{t}\}_{t\geq 0\text{ }}$ and $\{G_{t}\}_{t\in \mathbb{R} }$.

To prove the implication (iii)$\Rightarrow $(iv), we first observe
that if $\eta \in \partial \Delta $ is the Denjoy-Wolff point of
the group $ \{G_{t}\}_{t\in \mathbb{R} }$, then one can find an
automorphism $\psi $ of $\Delta $ such that $\psi (1)=\eta $. In
this case the family $\{\widetilde{G_{t}}\}_{t\in \mathbb{R} }$,
defined by
\begin{equation} \label{N4''}
\widetilde{G_{t}}=\psi ^{-1}\circ G_{t}\circ \psi ,
\end{equation}
has the Denjoy-Wolff point $\tau =1$.

Setting $\widetilde{\varphi }=\varphi \circ \psi $, we have by
(\ref{N2''}) and (\ref{N4''}) that $\widetilde{\varphi }(\Delta
)=\varphi (\Delta )=\Omega $ and
\begin{eqnarray*}
\widetilde{\varphi }(\widetilde{G_{t}}(z)) &=&\varphi (\psi
(\widetilde{G_{t} }(z))=\varphi (G_{t}(\psi (z)))= \\
&=&F_{t}(\varphi (\psi (z)))=F_{t}(\widetilde{\varphi
}(z))\text{.}
\end{eqnarray*}

So, in our setting, up to an automorphism we can assume that $\tau
=1$ is the Denjoy-Wolff point of the group $\{G_{t}\}_{t\in
\mathbb{R}}$ defined by (\ref{N1''}). Then it follows from
(\ref{N3''}) and Lindel\"{o}f's Theorem (see, for example,
\cite{G-L} and \cite{SD}) that
\begin{equation}  \label{N5''}
\angle \lim_{z\rightarrow 1}\varphi (z)=1\text{.}
\end{equation}

In addition, it is well known that the generator $g:\Delta \mapsto
\mathbb{C} $ of the group $\{G_{t}\}_{t\in \mathbb{R} }$ must be
of the form
\begin{equation} \label{N6''}
g(z)=a(z^{2}-1)+ib(z-1)^{2}
\end{equation}
with $a\geq 0$ and $b\in \mathbb{R} $. Moreover, $a=0$ if and only
if $\{G_{t}\}_{t\geq 0}$ is of parabolic type.

In the last case, for each $z\in \Delta $ we have
\begin{equation} \label{N7''}
\lim_{t\rightarrow \pm \infty }G_{t}(z)=1\text{.}
\end{equation}

If $a\neq 0$ $(a>0)$, then the point $\eta =-\dfrac{a+ib}{a-ib}$
is a repelling point of the group $\{G_{t}\}_{t\in \mathbb{R} }$,
i.e., $\lim_{t\rightarrow -\infty }G_{t}(z)=\eta $ for all $z\in
\Delta $ (or, which is the same, $g(\eta )=0$ and $g^{\prime
}(\eta )=-g^{\prime }(1)<0$). In this case, up to an automorphism
of $\Delta $ we can assume that $\eta =-1$ or, which is the same,
$b=0$.

In any case, differentiating (\ref{N2''}) at $t=0$ we obtain that
$\varphi $ must satisfy the differential equation
\begin{equation} \label{N8''}
\varphi ^{\prime }(z)\cdot g(z)=f(\varphi (z)),
\end{equation}
where $g$ is given by (\ref{N6''}).

Finally, we show that assertion (iv) implies (ii) as well as
(iii). Let $\varphi $ be a holomorphic solution of differential
equation (\ref{N8''}) with $\left\vert \varphi (z)\right\vert <1$,
where $g\in \Hol(\Delta , \mathbb{C} )$ is given by (\ref{N6''}).
We claim that $\varphi $ is univalent and is a conjugating
function for the semigroup $S=\{F_{t}\}_{t\geq 0\text{ }}$ with a
group $\{G_{t}\}_{t\in \mathbb{R}}$ of automorphisms of $\Delta $.

To this end, let us consider the function $h\in \Hol(\Delta ,
\mathbb{C} )$ defined by (\ref{D1}).

Recall that since $f(z)\neq 0$, $z\in \Delta $, the function $h$
is holomorphic in $\Delta $ and satisfies Abel's functional
equation
\begin{equation} \label{N10''}
h(F_{t}(z))=h(z)+t\text{, \ }z\in \Delta \text{, \ }t\geq 0\text{.}
\end{equation}

Let us define the function $k\in \Hol(\Delta , \mathbb{C} )$ by
the formula
\begin{equation} \label{N11''}
k(z)=h(\varphi (z)),\text{ \ \ }z\in \Delta \text{.}
\end{equation}

Then we have by (\ref{D1}) and (\ref{N11''})
\begin{equation*}
k^{\prime }(z)=h^{\prime }(\varphi (z))\cdot \varphi ^{\prime
}(z)=-\frac{ \varphi ^{\prime }(z)}{f(\varphi (z))}\text{.}
\end{equation*}

On the other hand, it follows from (\ref{N8''}) that
\begin{equation}  \label{N12''}
k^{\prime }(z)=-\frac{1}{g(z)}\text{.}
\end{equation}

Thus we get from (\ref{N6''})
\begin{equation}
k(z)=\left\{
\begin{array}{c}
\dfrac{1}{2a}\log \dfrac{1+\overline{\eta }z}{1-\overline{\eta
}z}+C,\text{ \ }a\neq 0 \\
\\
-\dfrac{ibz}{1-z}+C,\text{ \ \ }a=0\text{,}
\end{array}
\right.  \label{N13''}
\end{equation}
where $C=h(\varphi (0))$ and $\eta =-\frac{a+ib}{a-ib}$.

First, this formula shows that
\begin{equation} \label{N14''}
\varphi =h^{-1}\circ k
\end{equation}
must be univalent on $\Delta $.

Second, $h(\Delta )$ must contain at least a nonempty horisontal
strip or even a horizontal half-plane $k(\Delta )=h(\varphi
(\Delta ))$.

In fact, we will see below that $h(\Delta )$ cannot contain a
half-plane if $S=\{F_{t}\}_{t\geq 0\text{ }}$ is of hyperbolic
type.

In any case, the trajectory $\widetilde{F_{t}}(w)=h^{-1}(h(w)+t)$
is well defined for each $w\in \Omega $ and all $t\in \mathbb{R}$.

Hence, the family $\{\widetilde{F_{t}}\}_{t\in \mathbb{R} }$ forms
a group on $\Omega =\varphi (\Delta )=h^{-1}(k(\Delta ))$ which is
by (\ref{N10''}) an extension of the semigroup $S=\{F_{t}\}_{t\geq
0\text{ } } $ on $\Omega $.

Then setting again
\begin{equation*}
G_{t}(z)=\varphi ^{-1}(F_{t}(\varphi (z)),
\end{equation*}
we obtain a conjugation group $\{G_{t}\}_{t\geq 0}$ of
automorphisms on $ \Delta $, the generator of which is exactly
$g(z)=a(z^{2}-1)+ib(z-1)^{2}$.

So, $\Omega =\varphi (\Delta )$ is of $h$-type if and only if
$a\neq 0$ and is of $p$-type otherwise $(a=0, \, b\neq 0)$.

Now we show that if $w\in \Omega $ and
\begin{equation*}
\zeta =\lim_{t\rightarrow -\infty }F_{t}(w),
\end{equation*}
then $\zeta $ must lie on the boundary of $\Delta $.

Indeed, setting $w=\varphi (z)$ for some $z\in \Delta $ we have
\begin{equation*}
\lim_{t\rightarrow -\infty }F_{t}(\varphi (z))=\lim_{t\rightarrow -\infty
}h^{-1}(h(\varphi (z))+t)\text{.}
\end{equation*}
Since $h(\varphi (z))+t=k(z)+t$ tends to infinity as $t$ goes to
$-\infty $, we have by Lemma 2 in \cite{Shapiro}, p.162, that
$\zeta =\lim_{t\rightarrow -\infty }F_{t}(w)$ must lie on
$\partial \Delta $.

On the other hand,
\begin{equation*}
\lim_{t\rightarrow -\infty }F_{t}(\varphi (z))=\lim_{t\rightarrow -\infty
}\varphi (G_{t}(z))\text{.}
\end{equation*}
Thus, we have obtained that if $\eta \in \partial \Delta $ is a
fixed point of the group $\{G_{t}\}_{t\in \mathbb{R}}$, then
\begin{equation*}
\zeta =\lim_{z\rightarrow \eta }\varphi (z)=:\varphi (\eta )
\end{equation*}
exists and belongs to $\partial \Delta $. If, in particular, $\eta
=1$, then $\zeta $ must also be $1$ by (\ref{N5''}).
\end{proof39}

\subsection{The classes $\Sigma ^{\protect\alpha } \, (\mathcal{G}^{\protect\alpha
}) $ and corners of domains}

Now we will concentrate on geometric properties of a backward flow
invariant domain $\Omega $ of $p$-type. We show, in particular,
that under some smoothness conditions $\Omega $ has a corner of
opening $\pi / \alpha$ (at the point $z=1$) if and only if $f\in
\mathcal{G}^{\alpha }[1]$ (respectively, $h\in \Sigma ^{\alpha
}[1]$) for some $\alpha \in (1,2]$.

We have proved that if $\Omega $ is of $p$-type, then for all
$w\in \Omega $,
\begin{equation} \label{N15}
\lim_{t\rightarrow -\infty }F_{t}(w)=1\text{.}
\end{equation}
As we have already mentioned, it follows from Theorem 2.4 in
\cite{C-D-P} that if for at least one $w\in \Delta $ the
trajectory $\{F_{t}(w)\}_{t\in \mathbb{R} }$ is well defined and
condition (\ref{N15}) holds, then there is a horizontal half-plane
contained in $h(\Delta )$. Thus condition (\ref{N15}) is a
necessary and sufficient condition for existence of a BFID $\Omega
$ of $p$ -type.

So, we assume now that for some $z\in \Delta $, the straight line
$h(z)+t$ is contained in $h(\Delta )$ for all $t\in \mathbb{R} $
and
\begin{equation} \label{1a}
\lim_{t\rightarrow \pm \infty }h^{-1}(h(z)+t)=1\text{.}
\end{equation}

If, in addition, $h\in \Sigma ^{\alpha }[1]$, i.e.,
\begin{equation} \label{2a}
\lim_{z\rightarrow 1}(1-z)^{1+\alpha }h^{\prime }(z)=\mu \neq 0,\infty ,
\end{equation}
then we have
\begin{equation} \label{3a}
\lim_{z\rightarrow 1}\frac{f(z)}{(1-z)^{1+\alpha
}}=a=-\frac{1}{\mu }\neq 0,\infty \text{,}
\end{equation}
where, $\alpha \in (0,2]$ and $f(z)=-\dfrac{1}{h^{\prime }(z)}$ is the
generator of the semigroup $u(t,z):=h^{-1}(h(z)+t)$.

It follows now from the Cauchy problem
\begin{equation}
\left\{
\begin{array}{c}
\dfrac{du}{dt}+f(u)=0 \\
\\
u(t,z)=z\in \Delta ,
\end{array}
\right.  \label{4a}
\end{equation}
that
\begin{equation*}
\frac{du}{(1-u)^{1+\alpha }}=-\frac{f(u)dt}{(1-u)^{1+\alpha }}
\end{equation*}
or
\begin{equation*}
\frac{1}{\alpha }\cdot \frac{1}{t}\frac{1}{(1-u)^{\alpha
}}=-\frac{1}{t} \int\limits_{0}^{t}\frac{f(u)dt}{(1-u)^{1+\alpha
}}+\frac{1}{t(1-z)^{\alpha }}.
\end{equation*}
Hence,
\begin{equation*}
\lim_{t\rightarrow \pm \infty }\frac{1}{t(1-u(t,z))^{\alpha
}}=\frac{1}{\mu } \cdot \alpha
\end{equation*}
or
\begin{equation} \label{5a}
\lim_{t\rightarrow \pm \infty }t(1-u(t,z))^{\alpha }=\frac{\mu
}{\alpha } \neq 0,\infty \text{ .}
\end{equation}

Therefore,
\begin{equation} \label{5a'}
\lim_{t\rightarrow \infty }\arg (1-u(t,z))=\frac{1}{\alpha }\arg
\mu ,
\end{equation}
while
\begin{equation} \label{5b'}
\lim_{t\rightarrow -\infty }\arg (1-u(t,z))=\frac{1}{\alpha }(\arg
\mu \pm \pi )\text{.}
\end{equation}
Consequently,
\begin{equation}  \label{6a}
\left\vert \lim_{t\rightarrow \infty }\arg
(1-u(t,z))-\lim_{t\rightarrow -\infty }\arg (1-u(t,z))\right\vert
=\frac{\pi }{\alpha }\text{.}
\end{equation}

Since the trajectory $\{u(t,z)\}_{t\in \mathbb{R} }$ belongs to
$\Delta $, we see that $\alpha $ must be greater than or equal to
$1$.

This leads us to the following assertion.

\begin{proposition} \label{propA}
Let $f\in \mathcal{G}^{\alpha }[1]$ with
\begin{equation*}
\lim_{z\rightarrow 1}\frac{f(z)}{(1-z)^{1+\alpha }}=a\neq 0,\infty
\end{equation*}
and let $S=\{F_{t}\}_{t\geq 0}$ be the semigroup generated by $f$.
Assume that for a point $z\in \Delta $, the trajectory
$\{F_{t}(z)\}_{t\geq 0}$ can be extended to the whole real axis
and
\begin{equation*}
\lim_{t\rightarrow -\infty }(1-F_{t}(z))=1\text{.}
\end{equation*}
Then $\alpha \in \lbrack 1,2]$ and

(i) \ $\left\vert \arg (-a)\right\vert =\frac{\pi }{2}(2-\alpha )$;

(ii) the trajectory $\{F_{t}(z)\}_{t<0}$ converges tangentially to
the point $z=1$ as $t$ tends to $-\infty $. Moreover, if $\arg
(-a)\neq 0$, then
\begin{equation*}
\lim_{t\rightarrow -\infty }\arg (1-F_{t}(z))=(sign\arg
(-a))\frac{\pi }{2} \text{;}
\end{equation*}

(iii) there is a BFID $\Omega \subset \Delta $ of $p$-type which
has a corner of opening $\pi \gamma $ at the point $z=1$ with
$\gamma =\dfrac{1}{ \alpha }$.
\end{proposition}

\begin{proof} \label{proofN2}
Since $\{F_{t}(z)\}_{t\in \mathbb{R}}\subset \Delta $ we have
\begin{equation}  \label{n7a}
-\frac{\pi }{2}\leq \lim_{t\rightarrow -\infty }\arg
(1-F_{t}(z))\leq \frac{ \pi }{2}\text{.}
\end{equation}
Assume that $\arg (-a)(=-\arg \mu )<0$. Since $\alpha \leq 2$ we
have
\begin{equation*}
\frac{1}{\alpha }(\arg \mu +\pi )>\frac{\pi }{\alpha }\geq
\frac{\pi }{2} \text{.}
\end{equation*}

Hence the sign $(+)$ in (\ref{5b'}) is impossible and we have by
(\ref{5b'}) and (\ref{n7a}),
\begin{equation*}
\lim_{t\rightarrow -\infty }\arg (1-F_{t}(z))=\frac{1}{\alpha }(\arg \mu
-\pi )\geq -\frac{\pi }{2}\text{.}
\end{equation*}
This, in turn, implies that
\begin{equation*}
\arg \mu \geq \pi -\frac{\pi }{2}\alpha =\frac{\pi }{2}(2-\alpha )\text{.}
\end{equation*}

But it follows from Theorem \ref{teor3} and (\ref{3a}) that
\begin{equation*}
\arg \mu (=\left\vert \arg (-a)\right\vert )\leq \frac{\pi
}{2}(2-\alpha ) \text{.}
\end{equation*}
Thus
\begin{equation} \label{n7b}
\arg \mu =\frac{\pi }{2}(2-\alpha )\text{.}
\end{equation}
So we again obtain from (\ref{5b'}) that
\begin{equation*}
\lim_{t\rightarrow -\infty }\arg (1-F_{t}(z))=\frac{1}{\alpha
}\left[ \pi - \frac{\pi \alpha }{2}-\pi \right] =-\frac{\pi
}{2}\text{.}
\end{equation*}

In a similar way, we prove our assertion if $\arg (-a)(=-\arg \mu
)>0$.

If $\arg \mu =0$, then
\begin{equation*}
\lim_{t\rightarrow -\infty }\left\vert \arg (1-F_{t}(z)\right\vert
=\frac{ \pi }{\alpha }\geq \frac{\pi }{2}
\end{equation*}
and we get $\alpha =2$.
\end{proof}

Thus we have proved also the following assertions.

\begin{corollary} \label{colorA}
If the assumptions of Proposition \ref{propA} hold, then
\begin{equation*}
\lim_{t\rightarrow \infty }\arg (1-F_{t}(z))=-\mathrm{sgn}\arg
(-a)\cdot \frac{\pi (2-\alpha )}{2\alpha }\text{.}
\end{equation*}

Thus, the trajectory $\{F_{t}(z)\}_{t\geq 0}$ is tangential to the unit
circle at the point $z=1$ if and only if $\alpha =1$.

This trajectory is tangential to the real axis at the point $z=1$ if and
only if $\alpha =2$.
\end{corollary}

\begin{corollary} \label{corollary N}
Let $h\in \Sigma ^{\alpha }[1]$ with $\mu =\lim
(1-z)^{1+\alpha }h^{\prime }(z)$. The following assertions are equivalent:

(i) there is $z\in \Delta $ such that the point $h(z)+t\in
h(\Delta )$ for all $t\in \mathbb{R} $ and
\begin{equation*}
\lim_{t\rightarrow -\infty }h^{-1}(h(z)+t)=1\text{;}
\end{equation*}

(ii) $h(\Delta )\supset k(\Delta )$ (a half-plane), where
\begin{equation*}
k(z)=\frac{ibz}{1-z}+c,
\end{equation*}
for some $b,c\in \mathbb{R} $.

Moreover, in this case:

(a) $\alpha \in \lbrack 1,2]$;

(b) if $\alpha \in \lbrack 1,2)$, then there is no horizontal
half-plane $ H\subset h(\Delta )$ such that $H\cap k(\Delta
)=\emptyset$ and $b\cdot \arg \mu >0$.
\end{corollary}

\begin{proof} \label{proofD}
In fact, we only have to prove assertion (b). Indeed, if $ \alpha
\neq 2$, then there is only one maximal BFID $\Omega $ of
$p$-type, since its angle of opening at $z=1$ is $\dfrac{\pi
}{\alpha } > \dfrac{\pi }{2}$. Let now $\varphi =h^{-1}\circ k$ be
a conformal mapping on $\Omega $. We claim that the limit
\begin{equation*}
\angle \lim_{z\rightarrow 1}\frac{1-z}{[1-\varphi (z)]^{\alpha }}
\end{equation*}
exists finitely and is different from zero.

To prove this, we note that the group $\{G_{t}\}_{t\in \mathbb{R}
}$ is generated by $g$, $g(z)=ib(z-1)^{2}$, which belongs to the
class $\mathcal{G}^{1}[1]$ with $\lim_{z\rightarrow
1}\dfrac{g(z)}{(1-z)^{2}}=ib$. In addition, since $ F_{t}(\varphi
(z))=\varphi (G_{t}(z))$, $t\geq 0$, $z\in \Delta $, we have by
formula (\ref{5a}),
\begin{equation*}
\lim_{t\rightarrow \infty }\frac{1-G_{t}(z)}{[1-\varphi
(G_{t}(z))]^{\alpha } }=\lim_{t\rightarrow \infty
}\frac{t(1-G_{t}(z))}{t[1-F_{t}(\varphi (z))]^{\alpha
}}=-\frac{\alpha }{\mu ib}\text{.}
\end{equation*}

Since the function
\begin{equation*}
\frac{(1-z)^{1/\alpha }}{1-\varphi (z)}
\end{equation*}
does not admit negative real values, we get by the generalized
Lindel\"{o}f Theorem (see Theorem 2.20 in \cite{G-L}, p. 42) that
\begin{equation*}
\lim_{z\rightarrow 1}\frac{1-z}{[1-\varphi (z)]^{\alpha
}}=-\frac{\alpha }{ \mu ib}\neq 0,\infty \text{.}
\end{equation*}

Hence the limit
\begin{equation*}
m:=\angle \lim_{z\rightarrow 1}\frac{(1-z)^{1/\alpha }}{1-\varphi (z)}
\end{equation*}
is also finite and different from zero. Since the function
$\dfrac{1}{ 1-\varphi (z)}$ is of positive real part, it follows
from Lemma \ref{lem4} that
\begin{equation*}
\left\vert \arg m\right\vert \leq \frac{\pi }{2}(1-\frac{1}{\alpha })
\end{equation*}
or
\begin{eqnarray*}
\left\vert \arg m^{\alpha }\right\vert &=&\left\vert \arg \left(
-\frac{ \alpha }{\mu ib}\right) \right\vert = \\ &=&\left\vert
\frac{\pi }{2}-\arg (\mu b)\right\vert \leq \frac{\pi }{2} (\alpha
-1)\text{.}
\end{eqnarray*}
This implies that if $b>0$, then also $\arg \mu >0$ and
conversely: if $b<0$, then $\arg \mu >0$, and we are done.
\end{proof}

In some sense a converse assertion also holds.

\begin{theorem} \label{TheorB1}
Let $h\in \Sigma [1]$ and suppose that $ F_{t}(z)=h^{-1}(h(z)+t)$
can be extended for some $z\in \Delta $ to the whole real axis
with
\begin{equation*}
\lim_{t\rightarrow -\infty }F_{t}(z)=1\text{.}
\end{equation*}
In other words, there is a nonempty BFID $\Omega $ of $p$-type.
Suppose further that $\partial \Omega $ has a Dini-smooth corner
of opening $\pi \gamma $ at the point $z=1$. The following
assertions hold:

(i) $h\in \Sigma _{A}^{\alpha }[1]$, with $\alpha
=\dfrac{1}{\gamma }$ and $ \mu :=\angle \lim_{z\rightarrow
1}(1-z)^{1+\alpha }h^{\prime }(z)\neq 0,\infty $;

(ii) $\frac{1}{2}\leq \gamma \leq 1$, hence the angle of the
opening cannot be less than $\dfrac{\pi }{2}$ and $\left\vert \arg
\mu \right\vert \leq \dfrac{ \pi (2\gamma -1)}{2\gamma }$;

(iii) the angular limit of the inverse Visser-Ostrovskii quotient
\begin{equation*}
\angle \lim_{z\rightarrow 1}\frac{h(z)}{(z-1)h^{\prime }(z)}
\end{equation*}
exists and equals $\gamma $.
\end{theorem}

\begin{proof} \label{proofB1}
We have already shown that there is a Riemann conformal mapping
from $\Delta $ onto $\Omega $ such that $\varphi (1)=1$ and
$\varphi $ satisfies the differential equation
\begin{equation} \label{w1}
ib(1-z)^{2}\varphi ^{\prime }(z)=f(\varphi (z)),
\end{equation}
where $f=-\frac{1}{h}\in \mathcal{G}[1]$.

It follows from Theorem 3.9 in \cite{PC-92}, p. 52, that the
limits
\begin{equation} \label{w2}
\lim_{z\rightarrow 1}\frac{1-\varphi (z)}{(1-z)^{\gamma }}=m\neq 0,\infty
\end{equation}
and
\begin{equation} \label{w3}
\lim_{z\rightarrow 1}\frac{\varphi ^{\prime }(z)}{(1-z)^{\gamma
-1}}=n\neq 0,\infty
\end{equation}
exist.

In addition, as in the proof of this theorem, one can show that
the limit of the Visser-Ostrovskii quotient
\begin{equation*}
\frac{(1-z)\varphi ^{\prime }(z)}{1-\varphi (z)}
\end{equation*}
exists and equals $\gamma $:
\begin{equation} \label{w4}
\lim_{z\rightarrow 1}\frac{(1-z)\varphi ^{\prime }(z)}{1-\varphi
(z)}=\frac{n }{m}=\gamma \text{.}
\end{equation}

Taking some $\alpha >0$, consider the limit
\begin{equation} \label{w5}
a:=\lim_{z\rightarrow 1}\frac{f(\varphi (z))}{(1-\varphi
(z))^{1+\alpha }} \text{.}
\end{equation}

We have by (\ref{w1}),
\begin{eqnarray*}
a &=&ib\lim_{z\rightarrow 1}\frac{(1-z)^{2}\varphi ^{\prime
}(z)}{(1-\varphi (z))^{1+\alpha }}= \\ &=&ib\lim_{z\rightarrow
1}\frac{(1-z)\varphi ^{\prime }(z)}{1-\varphi (z)} \cdot
\lim_{z\rightarrow 1}\frac{1-z}{(1-\varphi (z))^{\alpha }}\text{.}
\end{eqnarray*}
Hence we infer from (\ref{w2}) and (\ref{w4}) that the limit in
(\ref{w5}) exists and is different from zero if and only if
$\alpha =\frac{1}{\gamma }$.

If this is indeed the case we have
\begin{equation} \label{w6}
-\frac{1}{\mu }=a=\frac{ib}{\alpha m^{\alpha }}=\frac{ib\gamma
}{m^{\alpha }} \text{.}
\end{equation}

Furthermore, if we represent $f$ by the Berkson-Porta
representation formula $f(z)=-(1-z)^{2}p(z)$, then we see that the
limit
\begin{equation}  \label{w7}
\lim_{\substack{ w\rightarrow 1  \\ w=\varphi (z)\in \Omega
}}\frac{p(w)}{ (1-w)^{\alpha -1}}=-a
\end{equation}
exists finitely and is different from zero. Since $\gamma \leq 1$,
we have $\alpha \geq 1$. In addition, it follows from a rigidity
theorem in \cite {E-L-R-S} that $\alpha \leq 2$. So, $0\leq \alpha
-1\leq 1$. For $\alpha \in \lbrack 1,2)$ we have
\begin{equation*}
\left\vert \arg \frac{p(z)}{(1-z)^{\alpha -1}}\right\vert =\left\vert \arg
p(z)-(\alpha -1)\arg (1-z)\right\vert \leq \frac{\pi }{2}\alpha <\pi \text{.}
\end{equation*}

Thus the function $\dfrac{p(z)}{(1-z)^{1-\alpha }}$ does not admit
negative real values. Then it follows from (\ref{w7}) and the
generalized Lindel\"{o}f Theorem (see Theorem 2.20 in \cite{G-L},
p. 42) that
\begin{equation}  \label{w8}
\lim_{z\rightarrow 1}\frac{p(z)}{(1-z)^{\alpha -1}}=-a
\end{equation}
whenever $1\leq \alpha <2$.

If $\alpha =2$, then one can show by using the Riesz-Herglotz
representation of $p$ that the function $\dfrac{p(z)}{1-z}$ is
bounded on each nontangential approach region in $\Delta $ with
vertex at $1$. So, one can again apply Lindel\"{o}f's Theorem to
show that the angular limit
\begin{equation*}
\angle \lim_{z\rightarrow 1}\frac{p(z)}{1-z}=-a
\end{equation*}
exists. Moreover, in this case $a$ is a negative real number.
Setting $\mu =- \frac{1}{a}$, we see that the limit
\begin{eqnarray*}
\angle \lim_{z\rightarrow 1}(1-z)^{1+\alpha }h^{\prime }(z) &=& \\
\left( =-\angle \lim_{z\rightarrow 1}\frac{(1-z)^{\alpha -1}}{p(z)}\right)
&=&\mu \neq 0,\infty \text{.}
\end{eqnarray*}

So $h\in \Sigma ^{\alpha }[1]$, and by Theorem 10,
\begin{equation*}
\left\vert \arg \mu \right\vert \leq \frac{\pi }{2}(2-\alpha )=\frac{\pi
(2\gamma -1)}{2\gamma }\text{.}
\end{equation*}

Assertions (i) and (ii) are proved.

To prove assertion (iii), we again consider equation\ (\ref{w1}),
rewriting it in the form
\begin{equation*}
h^{\prime }(\varphi (z))\cdot \varphi ^{\prime
}(z)=\frac{-1}{ib(1-z)^{2}} \text{.}
\end{equation*}

Integrating this equation, we get
\begin{equation} \label{w9}
h(\varphi (z))=k(z):=\frac{-1}{ib}\frac{z}{1-z}+h(\varphi (0))\text{.}
\end{equation}
Since $\varphi $ is an inner conjugating function for the
semigroup $S=\{F_{t}\}_{t\geq 0}$ and the group $\{G_{t}\}_{t\in
R}$, we have, for all $w\in \Omega $,
\begin{equation*}
F_{t}(w)=\varphi (G_{t}(\varphi ^{-1}(w)))
\end{equation*}
or
\begin{equation*}
\varphi ^{-1}(F_{t}(w))=G_{t}(\varphi ^{-1}(w))\text{.}
\end{equation*}

Consequently,
\begin{equation*}
\lim_{t\rightarrow \infty }\varphi ^{-1}(F_{t}(w))=1\text{.}
\end{equation*}

In other words, there is a trajectory converging to $1$ such that
the limit of $\varphi ^{-1}$ along this trajectory is also $1$.

Therefore, we have from (\ref{w2}) that
\begin{equation} \label{w11}
\lim_{t\rightarrow \infty }\frac{(1-F_{t}(w))^{\alpha }}{1-\varphi
^{-1}(F_{t}(w))}=\lim_{t\rightarrow \infty }\frac{(1-\varphi (G_{t}(\varphi
^{-1}(w)))^{\alpha }}{1-G_{t}(\varphi ^{-1}(w))}=m^{\alpha }\text{.}
\end{equation}

On the other hand, by direct calculations from (\ref{w9}) we get
\begin{eqnarray*}
1-\varphi ^{-1}(F_{t}(w)) &=&1-k^{-1}(h(F_{t}(w)))= \\
&=&\frac{1}{1-ib(h(F_{t}(w))-h(\varphi (0))}\text{.}
\end{eqnarray*}

Therefore, using (\ref{w11}), we get
\begin{eqnarray}
m^{\alpha } &=&\lim_{t\rightarrow \infty }(1-F_{t}(w))^{\alpha }\left[
1+ib(h(F_{t}(w))-h(\varphi (0)))\right] =  \label{w12} \\
&=&-ib\lim_{t\rightarrow \infty }(1-F_{t}(w))^{\alpha }h(F_{t}(w))\text{.}
\notag
\end{eqnarray}

Now we observe that since $h$ is univalent and $\alpha \geq 1$,
the function $(1-z)^{\alpha }h(z)$ is bounded on each
nontangential approach region $\{ z\in \mathbb{C}:\left\vert
1-z\right\vert \leq r(1-\left\vert z\right\vert )$, $r>1\}$.
Indeed, it follows from Koebe's inequality (see, for example,
\cite{DP} and \cite{PC-92}) that
\begin{eqnarray*}
\left\vert 1-z\right\vert ^{\alpha }\left\vert h(z)\right\vert &\leq
&\left\vert h^{\prime }(0)\right\vert \cdot \frac{\left\vert 1-z\right\vert
}{1-\left\vert z\right\vert }\cdot \left\vert 1-z\right\vert ^{\alpha -1}\leq
\\
&\leq &r\left\vert h^{\prime }(0)\right\vert \left\vert 1-z\right\vert
^{\alpha -1}\text{.}
\end{eqnarray*}
Therefore, again by the Lindel\"{o}f's theorem and (\ref{w12}), we
arrive at
\begin{equation*}
\angle \lim_{z\rightarrow 1}(1-z)^{\alpha }h(z)=-\frac{m^{\alpha
}}{ib}\text{ .}
\end{equation*}

Now it follows from (\ref{w6}) that
\begin{equation*}
\angle \lim_{z\rightarrow 1}\frac{h(z)}{(z-1)h^{\prime }(z)}=-\angle
\lim_{z\rightarrow 1}\frac{(1-z)^{\alpha }h(z)}{(1-z)^{1+\alpha }h^{\prime
}(z)}=\gamma \text{.}
\end{equation*}
Theorem 44 is proved.
\end{proof}

\vspace{12pt}

{\bf Acknowledgments.} Dmitry Khavinson's research was supported
by a grant from the National Science Foundation. Simeon Reich was
partially supported by the Fund for the Promotion of Research at
the Technion and by the Technion President's Research Fund. This
research is part of the European Science Foundation Networking
Programme HCAA.

\end{document}